\documentclass[11pt, reqno]{amsart}

\usepackage[utf8]{inputenc}
\usepackage[T1]{fontenc}
\usepackage[english]{babel}
\usepackage{amsmath, amsthm, amssymb, mathtools}
\usepackage[hmarginratio=1:1]{geometry}
\usepackage{enumitem}
\usepackage{tikz}
\usepackage[hyphens]{url}
\usepackage[hidelinks]{hyperref}

\newcommand{\R}{\mathbb{R}}
\newcommand{\N}{\mathbb{N}}
\newcommand{\Z}{\mathbb{Z}}
\newcommand{\C}{\mathbb{C}}
\newcommand{\dd}{\mathrm{d}}

\newtheorem{theorem}{Theorem}[section]
\newtheorem{lemma}[theorem]{Lemma}
\newtheorem{proposition}[theorem]{Proposition}
\newtheorem{corollary}[theorem]{Corollary}
\newtheorem{conjecture}[theorem]{Conjecture}
\theoremstyle{remark}
\newtheorem{remark}[theorem]{Remark}
\numberwithin{equation}{section}
\begin{document}

\title[Exact eigenvalue counts for the disk]{Exact eigenvalue counts for the disk}

\author{Matthias Baur}
\address{Institute of Analysis, Dynamics and Modeling, Department of Mathematics, University of Stuttgart, Pfaffenwaldring 57, 70569 Stuttgart, Germany}
\email{matthias.baur@mathematik.uni-stuttgart.de}


\begin{abstract}
We show that Weyl's law for the Dirichlet Laplacian on a disk $D_R$ of radius $R$ holds as an exact identity: the number of eigenvalues not exceeding $\lambda$ is equal to $$\frac{|D_R|}{4\pi}\lambda-\sum_{m=1}^\infty\{\nu_m(R\sqrt\lambda)\}$$ 
where
$\nu_1(x)>\nu_2(x)> ...$ are the zeros of the function $\nu\mapsto J_\nu(x)$, $J_\nu$ is the Bessel function of order $\nu$ and $\{\cdot\}$ denotes the fractional part.
P\'olya's conjecture for the disk is an immediate consequence. The identity is the zero-field limit of a similar exact formula for the
number of eigenvalues of the magnetic Dirichlet Laplacian with constant magnetic field on the disk up to any Landau level, which
is expressed through the zeros of Laguerre polynomials with respect to their parameter.
\end{abstract}

\maketitle

\section{Introduction and main results}

\subsection{P\'olya's conjecture and an exact Weyl formula for the disk}\label{sec:intro_weyl}

Let $\Omega\subset\R^2$ be an open set of finite measure and let $H_0^\Omega=-\Delta_\Omega^D$ denote the Dirichlet Laplacian on
$\Omega$. Its spectrum consists of eigenvalues $0<\lambda_1(\Omega)\le\lambda_2(\Omega)\le ...$, repeated according to
multiplicity, and we write
\[
N(\lambda,H_0^\Omega):=\#\{n\in\N:\lambda_n(\Omega)\le\lambda\},\qquad\lambda\in\R,
\]
for the eigenvalue counting function. 

Weyl's law \cite{Weyl} states that $N(\lambda,H_0^\Omega)=\frac{|\Omega|}{4\pi}\lambda(1+o(1))$ as $\lambda\to\infty$.
P\'olya conjectured that the leading term is in fact an upper bound for every $\lambda$, i.e.
\begin{equation}\label{eq:polya}
N(\lambda,H_0^\Omega)\le\frac{|\Omega|}{4\pi}\lambda,\qquad \lambda>0,
\end{equation}
or equivalently
\begin{equation}\label{eq:polya_lambda_n}
\lambda_n(\Omega)\ge \frac{4\pi n}{|\Omega|},  \qquad n\in \N,
\end{equation}
and proved this for domains which tile the plane \cite{Polya}. For bounded simply connected domains, a proof of the strict inequality has recently been announced in \cite{TangWang}. In general, the conjecture is still open. A particularly prominent
domain is the disk. It does not tile the plane, its eigenvalues are explicitly given in terms of zeros of Bessel functions, and
yet P\'olya's conjecture for the disk remained open for decades. It was finally proved by Filonov, Levitin, Polterovich and Sher
\cite{FLPS}, together with the conjecture for Dirichlet balls in any dimension, see also \cite{GMWZ,Li} for recent developments.
Their proof bounds the counting function by the number of lattice points under the graph of an explicit function, which rests on
uniform estimates for the zeros of Bessel functions, and then estimates this lattice point count.

In this paper we show that for the disk, Weyl's law can be turned into an exact identity. Throughout, $D_R$ denotes the open disk
of radius $R>0$ centred at the origin. Remarkably, the identity does not involve the zeros of the Bessel functions $x\mapsto J_\nu(x)$ of fixed
order $\nu$, but the zeros of $J_\nu(x)$ with respect to the order $\nu$. For fixed $x>0$, the function $\nu\mapsto J_\nu(x)$ has infinitely
many zeros. All of them are real and simple, and they are bounded from above, see Lemma~\ref{thm:bessel}, and we denote them by
$\nu_1(x)>\nu_2(x)>...$\,. They satisfy $\nu_m(x)>-m$ for all $m\in\N$, and $\nu_m(x)+m\to0$ as $m\to\infty$. Finally, for
$t\in\R$ let $\lfloor t\rfloor=\max\{k\in\Z:k\le t\}$ denote the largest integer not exceeding $t$ and $\{t\}=t-\lfloor t\rfloor$
the fractional part of $t$. The following identity holds.

\begin{theorem}[Exact Weyl formula for the disk]\label{thm:weyl}
Let $R>0$ and $\lambda>0$. Then
\begin{equation}\label{eq:weyl}
N(\lambda,H_0^{D_R})=\sum_{m=1}^{\infty}\Bigl(\bigl\lfloor\nu_m(R\sqrt\lambda)\bigr\rfloor+m\Bigr)
=\frac{|D_R|}{4\pi}\lambda-\sum_{m=1}^{\infty}\bigl\{\nu_m(R\sqrt\lambda)\bigr\},
\end{equation}
where only finitely many terms of the first series are different from zero and the second series converges.
\end{theorem}

In other words, for the disk the Weyl term $\frac{|D_R|}{4\pi}\lambda$ is exact up to a remainder which is explicit and has a
sign. Since $0<\nu_m(x)+m<1$, and hence $\{\nu_m(x)\}=\nu_m(x)+m>0$, for all sufficiently large $m$, the last series in
\eqref{eq:weyl} is strictly positive, and P\'olya's inequality follows at once.

\begin{corollary}[P\'olya's inequality for the disk]\label{cor:polya}
For every $R>0$,
\[
N(\lambda,H_0^{D_R})<\frac{|D_R|}{4\pi}\lambda,\qquad\lambda>0,
\]
or equivalently,
\[
\lambda_n(D_R)>\frac{4\pi n}{|D_R|}, \qquad n\in\N.
\]
\end{corollary}

The mechanism behind Theorem~\ref{thm:weyl} is a change of perspective. By separation of variables, the eigenvalues of
$H_0^{D_R}$ are the numbers $j_{|l|,k}^2/R^2$ with $l\in\Z$ and $k\in\N$, where $j_{\nu,k}$ denotes the $k$-th positive zero of
$J_\nu$. Instead of fixing the order and following the zeros in $x$, we consider the zero set of $(\nu,x)\mapsto J_\nu(x)$ in
$\R\times(0,\infty)$ as a whole. We show that it is the disjoint union of the graphs of strictly increasing functions
$x_m\colon(-m,\infty)\to(0,\infty)$, $m\in\N$, and that $\nu_m$ is the inverse function of $x_m$, see
Proposition~\ref{prop:besselcurves}. Since $J_{-l}=(-1)^lJ_l$ for $l\in\Z$, the eigenvalues of
the disk, counted with multiplicity, correspond exactly to the points $(\nu,x)$ on these curves whose first coordinate is an integer.
Hence $N(\lambda,H_0^{D_R})$ is the number of such points at height at most $x=R\sqrt\lambda$, and the $m$-th curve carries
$\lfloor\nu_m(x)\rfloor+m$ of them. This is the first identity in \eqref{eq:weyl}. The Weyl term then enters through the sum rule
\[
\sum_{m=1}^{\infty} (\nu_m(x)+m )=\frac{x^2}{4},\qquad x>0,
\]
see \eqref{eq:sumrule}, which we prove with the help of the argument principle and which appears to be new. In this picture, the
leading term in Weyl's law for the disk is a statement about the sum of the zeros $\nu_m(x)$, and the remainder records by how
much these zeros fail to be integers.

We stress that P\'olya's inequality for the disk, including its strict form, is due to \cite{FLPS}. What is new here is the
identity \eqref{eq:weyl} and the proof of the inequality which results from it. This proof uses only elementary properties of
Bessel functions and the argument principle. In particular, no asymptotic expansions or uniform bounds for zeros of Bessel
functions are needed, and all $\lambda>0$ are treated at once. The two approaches are nevertheless closely related. By the above,
the eigenvalue count for the disk is itself an exact lattice point count, and the lattice point count of \cite{FLPS} is obtained
from it by replacing the zeros $\nu_m(x)$ with explicit semiclassical approximations. We refer to Section~\ref{sec:extensions} for a discussion.

The identity \eqref{eq:weyl} also gives a new description of the remainder in Weyl's law. By the two-term Weyl asymptotics for the
disk \cite{KuznecovFedosov}, see also \cite{Ivrii}, it implies $\sum_{m=1}^{\infty}\{\nu_m(x)\}=\frac x2+o(x)$ as $x\to\infty$, and bounds on the remainder in the
two-term asymptotics are equivalent to bounds on the deviation of this sum of fractional parts from $\frac x2$.

\subsection{Constant magnetic fields and an exact count at the Landau levels}\label{sec:intro_landau}

Theorem~\ref{thm:weyl} is the zero-field limit of a formula for a problem with magnetic field, which is our second main result
and from which \eqref{eq:weyl} was originally derived. For $B\ge0$ we consider the magnetic Dirichlet Laplacian
$H_B^\Omega=(-i\nabla+A)^2$ with constant magnetic field of strength $B$, defined as the Friedrichs extension of the quadratic form
\[
a[u]=\int_\Omega|(-i\nabla+A)u|^2\,\dd x,\qquad u\in C_0^\infty(\Omega),\qquad A(x)=\frac B2\begin{pmatrix}-x_2\\ x_1\end{pmatrix}.
\]
The operator $H_B^\Omega$ is positive and has compact resolvent. In accordance with the notation above, we denote its eigenvalues,
repeated according to multiplicity, by $0<\lambda_1(\Omega,B)\le\lambda_2(\Omega,B)\le ...$, and its counting function by
$N(\lambda,H_B^\Omega):=\#\{n\in\N:\lambda_n(\Omega,B)\le\lambda\}$. If $B=0$, then $\lambda_n(\Omega,0) = \lambda_n(\Omega)$ are the eigenvalues of the usual Dirichlet Laplacian. The normalized magnetic flux through $\Omega$ is
$\phi=B|\Omega|/(2\pi)$, so that $\phi=BR^2/2$ for the disk $D_R$.

In the presence of a magnetic field, the analogue of \eqref{eq:polya} fails. A result of Erd\H{o}s, Loss and Vougalter \cite{ELV} implies that
\begin{align} \label{eq:ELV}
\lambda_n(\Omega,B)\ge\frac12\cdot\frac{4\pi n}{|\Omega|},   \qquad n\in \N,\ B\geq 0,
\end{align}
and the factor $1/2$ cannot be improved uniformly in $B$ \cite{BaurWeidl,FLW}. On the other hand, the spectrum of the Landau
Hamiltonian in the plane consists of the Landau levels $(2k-1)B$, $k\in\N$, each of which semiclassically carries $\phi$ states on a set of
measure $|\Omega|$. In accordance with this picture, it was shown in \cite[Theorem~2.2]{FLW} that for \emph{tiling} domains
$\Omega \subset \R^2$
\begin{equation}\label{eq:FLW}
\#\{n\in\N:\lambda_n(\Omega,B)<\lambda\}\ \le\ \#\{k\in\N:(2k-1)B<\lambda\} \cdot \phi,\qquad \lambda>0,\ B>0.
\end{equation}
The right-hand side resolves the Landau levels: it is a step function of $\lambda$ which increases by the flux $\phi$ at each level.
For $\lambda=(2j+1)B$ it equals $j\phi$, whereas the Weyl term $\frac{|\Omega|}{4\pi}\lambda$ equals $(j+\frac12)\phi$.

It is natural to ask whether \eqref{eq:FLW} holds beyond tiling domains. For the disk the answer is positive, again in the strong
form of an exact formula, now for the number of eigenvalues up to a Landau level. The role of the Bessel functions is taken by
Laguerre polynomials. For $j\in\N_0=\N\cup\{0\}$ and $\alpha\in\R$ let
\begin{equation}\label{eq:laguerre}
L_j^{(\alpha)}(z)=\sum_{i=0}^{j}\binom{j+\alpha}{j-i}\frac{(-z)^i}{i!}, \qquad z\in \R,
\end{equation}
be the Laguerre polynomial of degree $j$ with parameter $\alpha$. For fixed $z$, the right-hand side is also a polynomial of
degree $j$ in $\alpha$. We show in Lemma~\ref{cor:roots} that for $j\in\N$ and $z>0$ the polynomial
$\alpha\mapsto L_j^{(\alpha)}(z)$ has $j$ distinct real zeros $\alpha^{(j)}_1(z)>...>\alpha^{(j)}_j(z)$, and that they satisfy
$\alpha^{(j)}_m(z)>-m$ for $m=1,...,j$.

\begin{theorem}[Exact count at the Landau levels]\label{thm:landau}
Let $R>0$, $B>0$, $j\in\N$ and $\phi=BR^2/2$. Then
\begin{equation}\label{eq:landau}
N((2j+1)B,H_B^{D_R})=\sum_{m=1}^{j}\Bigl(\bigl\lfloor\alpha^{(j)}_m(\phi)\bigr\rfloor+m\Bigr)
=j\phi-\sum_{m=1}^{j}\bigl\{\alpha^{(j)}_m(\phi)\bigr\}.
\end{equation}
In particular,
\begin{equation}\label{eq:twosided}
j\phi-j<N((2j+1)B,H_B^{D_R})\le j\phi ,
\end{equation}
with equality on the right if and only if all $\alpha^{(j)}_m(\phi)$ are integers.
\end{theorem}

Theorem~\ref{thm:landau} says that the first $j$ Landau levels of the disk hold the semiclassically expected $j\phi$ states up to
a defect between $0$ and $j$ which is caused by the Dirichlet boundary condition and which is explicit. For example, for $j=1$ we
have $\alpha^{(1)}_1(\phi)=\phi-1$ and hence
\[
N(3B,H_B^{D_R})=\lfloor\phi\rfloor.
\]
As a consequence of Theorem~\ref{thm:landau}, the bound \eqref{eq:ELV} can be improved for the disk between any two consecutive
Landau levels, and \eqref{eq:FLW} holds for the disk.

\begin{corollary}\label{cor:flux}
Let $R>0$ and $B>0$.
\begin{enumerate}[label=\textup{(\roman*)}]
\item If $k\in\N$ and $(2k-1)B<\lambda_n(D_R,B)\le(2k+1)B$, then
\[
\lambda_n(D_R,B)>\Bigl(1-\frac{1}{2k}\Bigr)\frac{4\pi n}{|D_R|}.
\]
\item The disk $\Omega=D_R$ satisfies \eqref{eq:FLW}.
\end{enumerate}
\end{corollary}

Every eigenvalue of $H_B^{D_R}$ is covered by~(i), since $\lambda_1(D_R,B)>B$, and the factor $(1-\frac{1}{2k})$ is closer to one the more Landau levels lie below $\lambda_n(D_R,B)$. To the best of our knowledge, the disk is thus the first non-tiling domain for which \eqref{eq:FLW} is known.

The two theorems are instances of the same picture. In the proof of Theorem~\ref{thm:landau}, the zero set of
$(\alpha,z)\mapsto L_j^{(\alpha)}(z)$ takes the role of the zero set of $(\nu,x)\mapsto J_\nu(x)$. It consists of $j$ increasing
curves, the eigenvalues up to the Landau level $(2j+1)B$ correspond to the points on these curves with integer $\alpha$ at height
at most $\phi$, and the term $j\phi$ in \eqref{eq:landau} stems from Vieta's formula
$\sum_{m=1}^{j} (\alpha^{(j)}_m(\phi)+m )=j\phi$ for the sum of the zeros. If one lets $j\to\infty$ along
$\phi=x^2/(4j)$ with fixed $x=R\sqrt\lambda>0$, then $(2j+1)B\to\lambda$ and $j\phi=x^2/4$, suitably rescaled Laguerre polynomials
converge to Bessel functions, $\alpha^{(j)}_m(\phi)\to\nu_m(x)$, and \eqref{eq:landau} turns into \eqref{eq:weyl}. In this sense,
P\'olya's inequality for the disk is the zero-field limit of the upper bound \eqref{eq:twosided} at the Landau levels. Its
non-strict version already follows from Corollary~\ref{cor:flux}\,(i) by letting $B\to0$, see Remark~\ref{rem:Bto0}. However, our proof of
Theorem~\ref{thm:weyl} does not pass to the limit in the counting functions. We give two proofs of the required properties of the
zeros $\nu_m(x)$: one in which reality and monotonicity are inherited from the zeros of Laguerre polynomials, and an independent
one which realizes the numbers $\nu_m(x)$ as the eigenvalues of a Jacobi operator, see Remark~\ref{rem:jacobi}.

\subsection{Further results and open questions}\label{sec:intro_further}

The description of the zero set of $(\nu,x)\mapsto J_\nu(x)$ applies to other spectral problems whose eigenvalues are given by
zeros of Bessel functions. In Section~\ref{sec:extensions} we derive exact counting formulas in terms of the zeros $\nu_m(x)$ for
the disk with an Aharonov-Bohm flux and for the Dirichlet Laplacian on balls in $\R^d$, and we explain how the lattice point
count of \cite{FLPS} fits into this picture.

For Neumann boundary conditions we expect analogous exact formulas, with and without magnetic field. In the latter case, the zeros of $\nu\mapsto J_\nu(x)$ are replaced by the zeros of $\nu\mapsto J_\nu'(x)$ and the remainder has the opposite sign. We state them in Section~\ref{sec:neumann}. Since the corresponding zero curves are no longer monotone, the proofs require additional arguments, which we plan to give in a follow-up article.

Let us finally point out the limitations of our approach. It relies on the separation of variables and is therefore tied to the
disk and to closely related rotationally symmetric problems. Moreover, although our counting formula for balls is exact, we do not know whether it leads to a similarly direct proof of P\'olya's conjecture for balls in dimension $d\ge3$, which has also already been treated in \cite{FLPS}. This is a direction of future research.

\subsection{Organization of the paper and notation}

Sections~\ref{sec:laguerre} and~\ref{sec:landau} are concerned with the magnetic Dirichlet Laplacian on the disk. It is instructive to start with the magnetic case, since the proofs involve only polynomials and are thus a bit simpler than in the non-magnetic case. In
Section~\ref{sec:laguerre} we describe the zero set of $(\alpha,z)\mapsto L_j^{(\alpha)}(z)$, and Section~\ref{sec:landau} contains
the reduction of the eigenvalue count to zeros of Laguerre polynomials and the proofs of Theorem~\ref{thm:landau} and
Corollary~\ref{cor:flux}. Sections~\ref{sec:bessel} and~\ref{sec:weyl} are then concerned with the non-magnetic case. In
Section~\ref{sec:bessel} we describe the zero set of $(\nu,x)\mapsto J_\nu(x)$ and prove the sum rule \eqref{eq:sumrule}, parallel
to Section~\ref{sec:laguerre}, and in Section~\ref{sec:weyl} we prove Theorem~\ref{thm:weyl} and Corollary~\ref{cor:polya}.
Section~\ref{sec:extensions} contains the further results mentioned above. A reader who is only interested in
Theorem~\ref{thm:weyl} can read Section~\ref{sec:bessel} together with the alternative proofs in Remark~\ref{rem:jacobi}, followed by
Section~\ref{sec:weyl}, essentially independently of Sections~\ref{sec:laguerre} and~\ref{sec:landau}.

\section{Zeros of Laguerre polynomials as functions of the parameter}\label{sec:laguerre}

Let $(c)_k=c(c+1)\cdots(c+k-1)$, $(c)_0=1$, for $c\in\R$, $k\in\N_0$ denote the Pochhammer symbol. For $j\in\N_0$ and $(\alpha,z)\in\R^2$ we put
\begin{equation}\label{eq:Pj}
P_j(\alpha,z):=j! \, L_j^{(\alpha)}(z)=\sum_{i=0}^{j}\binom{j}{i}(\alpha+i+1)_{j-i}\,(-z)^i ,
\end{equation}
so that $P_0\equiv1$. The second expression follows from \eqref{eq:laguerre} since $\binom{j+\alpha}{j-i}=(\alpha+i+1)_{j-i}/(j-i)!$.

\begin{lemma}\label{lem:Pj}
Let $j\in\N$.
\begin{enumerate}[label=\textup{(\roman*)}]
\item For fixed $z\in\R$, $\alpha\mapsto P_j(\alpha,z)$ is a polynomial in $\alpha$ with
\begin{align*}
    P_j(\alpha,z) = \alpha^j + \bigg( \frac{j(j+1)}{2}-jz \bigg) \alpha^{j-1} + ... \,.
\end{align*}
\item For every $\alpha\in\R$ the function $y=P_j(\alpha,\cdot)$ satisfies $zy''+(\alpha+1-z)y'+jy=0$. Consequently, all zeros of
$P_j(\alpha,\cdot)$ in $(0,\infty)$ are simple.
\item $P_j(\alpha,0)=(\alpha+1)_j$, which vanishes if and only if $\alpha\in\{-1,-2,...,-j\}$.
\item If $\alpha>-1$, then $P_j(\alpha,\cdot)$ has exactly $j$ zeros in $(0,\infty)$. If $\alpha<-j$, it has no zeros in $(0,\infty)$.
\item For $k\in\{1,...,j\}$ and $z\in\R$,
\begin{equation}\label{eq:negparam}
P_j(-k,z)=(-z)^k\,P_{j-k}(k,z).
\end{equation}
In particular, $P_j(-k,\cdot)$ has exactly $j-k$ zeros in $(0,\infty)$, namely those of $L_{j-k}^{(k)}$.
\end{enumerate}
\end{lemma}

\begin{proof}
(i) This is a matter of collecting the coefficients of $\alpha^j$ and $\alpha^{j-1}$ in \eqref{eq:Pj}. The term of \eqref{eq:Pj} with $i=0$ is $(\alpha+1)_j=\alpha^j+\frac{j(j+1)}{2}\alpha^{j-1}+...$, the term with $i=1$ is $-jz(\alpha+2)_{j-1}=-jz\alpha^{j-1}+...$, and the terms with $i\ge2$ have degree at most $j-2$ in $\alpha$.

(ii) For $\alpha>-1$ this is the differential equation of the Laguerre polynomials, see
\cite[(5.1.2)]{Szego} or \cite[Table~18.8.1]{DLMF}. For fixed $z$, the expression
$zy''+(\alpha+1-z)y'+jy$ is a polynomial in $\alpha$ by \eqref{eq:Pj}, so the equation holds
for all $\alpha\in\R$. The differential equation is regular on $(0,\infty)$. Suppose a zero $z_0>0$ of
$P_j(\alpha,\cdot)$ is not simple. Then $y(z_0)=y'(z_0)=0$ and $y=P_j(\alpha,\cdot)\equiv0$ by uniqueness, which is a contradiction.

(iii) follows immediately from \eqref{eq:Pj}.

(iv) For $\alpha>-1$ the polynomials $L^{(\alpha)}_j$ are orthogonal with respect to the weight $z^\alpha e^{-z}$ on $(0,\infty)$,
so $L_j^{(\alpha)}$ has $j$ simple zeros in $(0,\infty)$ \cite[Theorem~3.3.1]{Szego}. If $\alpha<-j$, then all factors of
$(\alpha+i+1)_{j-i}=\prod_{t=i+1}^{j}(\alpha+t)$ are negative, so the sign of $(\alpha+i+1)_{j-i}$ is $(-1)^{j-i}$. Thus, for $z>0$ every summand in \eqref{eq:Pj} is $\neq 0$ and has sign $(-1)^j$. Hence, $P_j(\alpha,\cdot)$ has no zeros in $(0,\infty)$.

(v) For $\alpha=-k$ we have $(i+1-k)_{j-i}=\prod_{t=i+1}^{j}(t-k)$, which vanishes for $i<k$ and equals $(j-k)!/(i-k)!$ for
$i\ge k$. Hence, with $i=i'+k$,
\begin{align*}
P_j(-k,z)&=\sum_{i=k}^{j}\binom{j}{i} \frac{(j-k)!}{(i-k)!}\,(-z)^i \\
&=(-z)^k\sum_{i'=0}^{j-k}\frac{j!\,(j-k)!}{(i'+k)!\,(j-k-i')!\,i'!}\,(-z)^{i'}\\
&=(-z)^k\sum_{i'=0}^{j-k}\binom{j-k}{i'}(k+i'+1)_{j-k-i'}\,(-z)^{i'},
\end{align*}
because $(k+i'+1)_{j-k-i'}=j!/(i'+k)!$. This is \eqref{eq:negparam}. The last assertion follows from (iv), applied to $P_{j-k}(k,\cdot)$ if $k<j$, and from $P_0\equiv1$ if $k=j$.
\end{proof}

The identity \eqref{eq:negparam} is classical, see \cite[\S5.2]{Szego}. The next proposition is the key to all counting arguments below.

\begin{proposition}\label{prop:curves}
Let $j\in\N$ and
$$Z:=\{(\alpha,z)\in\R\times(0,\infty):P_j(\alpha,z)=0\}.$$
There exist continuous, strictly increasing bijections
\begin{align*}
    \zeta_m\colon(-m,\infty)&\to(0,\infty),\qquad m=1,...,j, \\
\alpha &\mapsto \zeta_m(\alpha),
\end{align*}
with the following properties:
\begin{enumerate}[label=\textup{(\roman*)}]
\item $ Z=\bigcup_{m=1}^{j}\bigl\{(\alpha,\zeta_m(\alpha)):\alpha>-m\bigr\}$.
\item If $1\le m<m'\le j$, then $\zeta_m(\alpha)<\zeta_{m'}(\alpha)$ for $\alpha>-m$.
\end{enumerate}
In other words, the zero set $Z$ consists precisely of $j$ non-intersecting, strictly increasing graphs.
\end{proposition}

\begin{proof} 

By Lemma~\ref{lem:Pj}\,(ii), $\partial_zP_j\neq0$ on $Z$. By the implicit function theorem, around every point of $Z$, $Z$ is locally the graph of a continuous function. Uniqueness implies that if $I$ is an interval and $\zeta,\tilde{\zeta}\colon I\to(0,\infty)$ are continuous functions whose graphs are contained in $Z$ and which coincide at one point, then $\zeta=\tilde{\zeta}$ on $I$. This uniqueness principle will be used throughout the proof.

\textit{Step 1: Construct $\zeta_m$.} Let us construct the functions $\zeta_m$. For $\alpha>-1$, $P_j(\alpha,\cdot)$ has exactly $j$ zeros in $(0,\infty)$, see Lemma~\ref{lem:Pj}\,(iv), so let $\zeta_1(\alpha)< ... <\zeta_j(\alpha)$ denote these zeros. Each $\zeta_m$ is continuous on $(-1,\infty)$. Indeed, for
$\alpha$ near an $\alpha_0>-1$, the local graphs through the points $(\alpha_0,\zeta_m(\alpha_0))$, $m=1,...,j$, provide locally the $j$ distinct zeros of $P_j(\alpha,\cdot)$ in the same order, so each $\zeta_m$ is continuous at $\alpha_0$. By patching together local graphs, we extend each $\zeta_m$ to the left as far as possible. Let $(a_m,\infty)$ be the largest interval to which $\zeta_m$ can be extended as a continuous function with graph in $Z$. By the uniqueness principle, the extension is unique and we denote it again by $\zeta_m$. Since $P_j(\alpha,\cdot)$ has no positive zeros for $\alpha<-j$ by Lemma~\ref{lem:Pj}\,(iv), we have $a_m\ge-j$ for every $m=1,...,j$.

\textit{Step 2: Limit of $\zeta_m$ as $\alpha\searrow a_m$.} We claim that 
\begin{align} \label{eq:prop2_amlimit}
  \zeta_m(\alpha)\to0 , \qquad \alpha\searrow a_m.
\end{align}
For $\alpha$ in a bounded set, the zeros of $P_j(\alpha,\cdot)$ are bounded, because the leading coefficient in $z$ is $(-1)^j$ and the other coefficients are bounded. In particular, $\zeta_m$ is bounded on $(a_m,a_m+1]$. If $\zeta_m(\alpha)\not\to 0$ as $\alpha\searrow a_m$, there would be $\varepsilon>0$ and a sequence $\alpha_i\searrow a_m$ with $\zeta_m(\alpha_i)\ge\varepsilon$, and, by compactness, after passing to a subsequence we may assume that $\zeta_m(\alpha_i)\to z^*$ for some $z^*\ge\varepsilon>0$. But  then, by continuity of $P_j$ with respect to both variables $(\alpha,z)$, we must have $(a_m,z^*)\in Z$, and by the implicit function theorem, there is a neighbourhood of $(a_m,z^*)$ and a local graph which extends $\zeta_m$ to the left of $a_m$. This contradicts the maximality of $(a_m,\infty)$ and proves the claim. Again by continuity, it follows that $P_j(a_m,0)=0$, so that
\begin{equation}\label{eq:am}
a_m\in\{-1,-2,...,-j\}
\end{equation}
by Lemma~\ref{lem:Pj}\,(iii).

\textit{Step 3: Limit of $\zeta_m$ as $\alpha\to \infty$.} We next claim that
\begin{align} \label{eq:prop2_inflimit}
  \zeta_m(\alpha) \to \infty , \qquad \alpha \to \infty.
\end{align}
Fix $c>0$. For $\alpha>-1$ and $z\in[0,c]$ we have by \eqref{eq:Pj}
\[
P_j(\alpha,z)\ge(\alpha+1)_j-\sum_{i=1}^{j}\binom ji(\alpha+i+1)_{j-i}c^i ,
\]
and the right-hand side is a polynomial in $\alpha$ with leading term $\alpha^j$. Hence $P_j(\alpha,\cdot)>0$ on $[0,c]$ for all
sufficiently large $\alpha$, so $\zeta_m(\alpha)>c$ for all
sufficiently large $\alpha$. Therefore, $\zeta_m(\alpha)\to\infty$ as $\alpha\to\infty$. By \eqref{eq:prop2_amlimit}, \eqref{eq:prop2_inflimit} and continuity, we have shown that each $\zeta_m$ maps $(a_m,\infty)$ onto $(0,\infty)$.

\emph{Step 4: Strict monotonicity of $\zeta_m$.} The graphs of $\zeta_1,...,\zeta_j$ are pairwise disjoint: they are disjoint over $(-1,\infty)$, and two of the functions cannot coincide at a point of their common domain by the uniqueness principle. Consequently, by continuity,
\begin{equation}\label{eq:order}
\zeta_m(\alpha)<\zeta_{m'}(\alpha)\qquad\text{for }m<m'\text{ and }\alpha>\max\{a_m,a_{m'}\}.
\end{equation}
Now let $z>0$ be arbitrary but fixed. Since $\zeta_m$ maps $(a_m,\infty)$ continuously onto $(0,\infty)$, there is $\alpha^{(m)}>a_m$ with $\zeta_m(\alpha^{(m)})=z$ for each $m=1,..., j$, and the numbers
$\alpha^{(1)},...,\alpha^{(j)}$ are pairwise distinct because the graphs are disjoint. These are $j$ distinct zeros of the polynomial $\alpha\mapsto P_j(\alpha,z)$, which has degree $j$ by Lemma~\ref{lem:Pj}\,(i). Hence there are no further zeros and this shows that each $\zeta_m$ attains the value $z$ exactly once. It follows that each $\zeta_m$ is injective, hence strictly monotone, hence strictly increasing by \eqref{eq:prop2_amlimit} and \eqref{eq:prop2_inflimit}. Moreover, every point of $Z$ lies on one of the graphs: if $(\alpha_0,z)\in Z$, then
$\alpha_0$ is a zero of $\alpha\mapsto P_j(\alpha,z)$, hence $\alpha_0=\alpha^{(m)}$ for some $m$,
which means $\alpha_0>a_m$ and $\zeta_m(\alpha_0)=z$. Since the graphs are contained in $Z$ by
construction, we obtain
\begin{equation}\label{eq:Zunion}
Z=\bigcup_{m=1}^{j}\bigl\{(\alpha,\zeta_m(\alpha)):\alpha>a_m\bigr\}.
\end{equation}

\emph{Step 5: $a_m=-m$.} Suppose $m<m'$ but $a_m<a_{m'}\leq -1$. Then \eqref{eq:order} gives $0<\zeta_m(\alpha)<\zeta_{m'}(\alpha)$ for
$\alpha>a_{m'}$, and letting $\alpha\searrow a_{m'}$ we obtain $\zeta_m(a_{m'})=0$, a contradiction to $\zeta_m(a_{m'})>0$, which holds since $a_m<a_{m'}$. Hence $a_1\ge a_2\ge ... \ge a_j$. Let $k\in\{0,1,...,j-1\}$. By \eqref{eq:Zunion}, the number of zeros of $P_j(-k,\cdot)$ in $(0,\infty)$ equals $\#\{m:a_m<-k\}$,
and by Lemma~\ref{lem:Pj}\,(iv) and (v) this number is $j-k$. Since the $a_m$ are non-increasing, this means $a_m<-k$ exactly for
$m\ge k+1$. Taking $k=m-1$ gives $a_m<-(m-1)$, and taking $k=m$ gives $a_m\ge-m$ (if $m\le j-1$). For $m=j$ we have $a_j\ge-j$ by Step~1. In view of \eqref{eq:am} we conclude $a_m=-m$ for every $m=1,..., j$. 

Finally, (i) is \eqref{eq:Zunion} and (ii) is \eqref{eq:order}.
\end{proof}

\begin{lemma}\label{cor:roots}
Let $j\in\N$ and $z>0$, and let $\alpha^{(j)}_m:=\zeta_m^{-1}\colon(0,\infty)\to(-m,\infty)$, $m=1,...,j$, be the inverse functions
of the functions $\zeta_m$ from Proposition~\ref{prop:curves}.
\begin{enumerate}[label=\textup{(\roman*)}]
\item The zeros of the polynomial $\alpha\mapsto P_j(\alpha,z)$ are exactly $\alpha^{(j)}_1(z),...,\alpha^{(j)}_j(z)$.
They are real and simple, and
\[
\alpha^{(j)}_1(z)>\alpha^{(j)}_2(z)>...>\alpha^{(j)}_j(z),\qquad\alpha^{(j)}_m(z)>-m .
\]
Each $\alpha^{(j)}_m$ is continuous and strictly increasing.
\item For every $\alpha\in\R$,
\[
\#\{x\in(0,z]:P_j(\alpha,x)=0\}=\#\{m\in\{1,...,j\}:-m<\alpha\le\alpha^{(j)}_m(z)\}.
\]
\item We have
\begin{align} \label{eq:vieta_alpha}
    \sum_{m=1}^{j}\alpha^{(j)}_m(z)=jz-\frac{j(j+1)}{2} .
\end{align}
\end{enumerate}
\end{lemma}

\begin{proof}
(i) By Proposition~\ref{prop:curves}\,(i), $P_j(\alpha,z)=0$ if and only if $z=\zeta_m(\alpha)$ for some $m$ with $\alpha>-m$, or equivalently, $\alpha=\alpha^{(j)}_m(z)$. If $m<m'$, then $\alpha^{(j)}_m(z)>-m>-m'$ lies in the domain of $\zeta_{m'}$, and
$\zeta_{m'}(\alpha^{(j)}_m(z))>\zeta_m(\alpha^{(j)}_m(z))=z=\zeta_{m'}(\alpha^{(j)}_{m'}(z))$ by Proposition~\ref{prop:curves}\,(ii). Since
$\zeta_{m'}$ is increasing, this implies $\alpha^{(j)}_m(z)>\alpha^{(j)}_{m'}(z)$ and thus the ordering of the zeros $\alpha^{(j)}_m(z)$. Furthermore, a polynomial of degree $j$ with $j$ distinct zeros has only simple zeros, so all zeros $\alpha^{(j)}_m(z)$ are simple.

(ii) By Proposition~\ref{prop:curves}\,(i), the zeros of $P_j(\alpha,\cdot)$ in $(0,\infty)$ are the pairwise distinct numbers
$\zeta_m(\alpha)$ with $-m<\alpha$, and $\zeta_m(\alpha)\le z$ is equivalent to $\alpha\le\alpha^{(j)}_m(z)$ by monotonicity.

(iii) follows from (i), Vieta's formula and Lemma~\ref{lem:Pj}\,(i).
\end{proof}

\section{The exact count at the Landau levels}\label{sec:landau}

Let 
$$M(a,b,z)=\sum_{k =0}^\infty\frac{(a)_k}{(b)_k\,k!}z^k$$
be Kummer's confluent hypergeometric function \cite[\S13.2]{DLMF}. We recall the explicit solution of the eigenvalue problem on
the disk from \cite[Section~2.1]{BaurWeidl}. Let $R>0$, $B>0$ and $\phi=BR^2/2$. For every $l\in\Z$ the real solutions $\lambda$ of
the equation
\begin{equation}\label{eq:kummer}
M\biggl(\frac12\Bigl(l+|l|+1-\frac{\lambda}{B}\Bigr),|l|+1,\phi\biggr)=0
\end{equation}
form an increasing sequence $\lambda_{1,l}(B)<\lambda_{2,l}(B)<...$. The numbers $\lambda_{m,l}(B)$, $m\in\N$, $l\in\Z$, are
precisely the eigenvalues of $H_B^{D_R}$, counted with multiplicity, where $l$ is the angular momentum of the corresponding
eigenfunction. Consequently,
\begin{equation}\label{eq:Nbranches}
N(\lambda,H_B^{D_R})=\#\{(m,l)\in\N\times\Z:\lambda_{m,l}(B)\le\lambda\},\qquad\lambda\in\R .
\end{equation}
We need the following facts about the eigenvalue branches $B\mapsto\lambda_{m,l}(B)$. For all $m\in\N$, $l\in\Z$ and any $B>0$,
\begin{equation}\label{eq:BWlower}
\lambda_{m,l}(B)>(l+|l|+1+2(m-1))B
\end{equation}
by \cite[Theorem~2.1]{BaurWeidl}, and the function $B\mapsto\lambda_{m,l}(B)/B$ is strictly decreasing on $(0,\infty)$, see
\cite[Section~3]{BaurWeidl}. Furthermore, for $k\in\N$ and $l\in\Z$ a comparison of the power series of $M$ with \eqref{eq:Pj} gives
\begin{equation}\label{eq:MP}
(|l|+1)_k\,M(-k,|l|+1,z)=P_k(|l|,z).
\end{equation}
By Lemma~\ref{lem:Pj}\,(iv), the polynomial $z\mapsto M(-k,|l|+1,z)$ therefore has exactly $k$ zeros in $(0,\infty)$, which we
denote by $z_1(-k,|l|+1)< ... <z_k(-k,|l|+1)$. The intersections of the eigenvalue branches with the lines
$\lambda=(l+|l|+1+2k)B$ are characterized by \cite[Corollary~4.1]{BaurWeidl}. If $m,k\in\N$ with $k\ge m$ and $l\in\Z$, then
\begin{equation}\label{eq:BWcor}
\lambda_{m,l}(B)=(l+|l|+1+2k)B\qquad\text{if and only if}\qquad \phi = z_m(-k,|l|+1).
\end{equation}

For the statement of the next lemma, let $t_+=\max\{t,0\}$ for $t\in\R$. 

\begin{lemma}\label{lem:sturm}
Let $R>0$, $B>0$, $\phi=BR^2/2$, $j\in\N_0$, $l\in\Z$ and $m\in\N$, and put $k:=j-l_+$. Then
\[
\lambda_{m,l}(B)\le(2j+1)B\qquad\text{if and only if}\qquad k\ge m\ \text{ and } \phi \geq z_m(-k,|l|+1) .
\]
\end{lemma}

\begin{proof}
Since $l+|l|=2l_+$, we have $(2j+1)B=(l+|l|+1+2k)B$. If $k\le m-1$, then \eqref{eq:BWlower} gives
$\lambda_{m,l}(B)>(l+|l|+1+2(m-1))B\ge(2j+1)B$. Now let $k\ge m$ and put $B_*:=2z_m(-k,|l|+1)/R^2$. By \eqref{eq:BWcor},
$\lambda_{m,l}(B_*)/B_*=2j+1$. Since $B\mapsto\lambda_{m,l}(B)/B$ is strictly decreasing, $\lambda_{m,l}(B)/B\le2j+1$ holds if and
only if $B\ge B_*$, that is, if and only if $\phi\ge z_m(-k,|l|+1)$.
\end{proof}

\begin{proof}[Proof of Theorem~\ref{thm:landau}]
Let $j\in\N_0$. If $l\ge j$, then $j-l_+\le0$ and Lemma~\ref{lem:sturm} shows that $\lambda_{m,l}(B)>(2j+1)B$ for all $m\in\N$.
Hence, by \eqref{eq:Nbranches} and Lemma~\ref{lem:sturm},
\begin{equation}\label{eq:Nsplit}
N((2j+1)B,H_B^{D_R})=\sum_{l\in\Z,\ l<j}\#\bigl\{m\in\{1,...,j-l_+\}:z_m(-(j-l_+),|l|+1)\le\phi\bigr\}.
\end{equation}
In particular, when $j=0$, then $N(B,H_B^{D_R})=0$.

Now let $j\in\N$ and $l<j$, and put $k=j-l_+\in\{1,...,j\}$. By \eqref{eq:MP}, the summand in \eqref{eq:Nsplit} is the number of
zeros of $P_k(|l|,\cdot)$ in $(0,\phi]$. For $l\le0$ we have $P_k(|l|,z)=P_j(-l,z)$. For $1\le l\le j-1$ we have
$P_k(|l|,z)=P_{j-l}(l,z)$, whose zeros in $(0,\infty)$ coincide with those of $P_j(-l,\cdot)$ by \eqref{eq:negparam}. In both cases,
the summand is the number of zeros of $P_j(-l,\cdot)$ in $(0,\phi]$. Writing $\alpha=-l$, we obtain
\begin{equation}\label{eq:Nall}
N((2j+1)B,H_B^{D_R})=\sum_{\alpha\in\Z,\ \alpha>-j}\#\{x\in(0,\phi]:P_j(\alpha,x)=0\}.
\end{equation}
By Lemma~\ref{cor:roots}\,(ii) and interchange of the order of counting,
\begin{align*}
N((2j+1)B,H_B^{D_R})&=\sum_{\alpha\in\Z,\ \alpha>-j}\# \{m\in\{1,...,j\}:-m<\alpha\le\alpha^{(j)}_m(\phi) \} \\
&=\sum_{m=1}^{j}\# \{\alpha\in\Z:-m<\alpha\le\alpha^{(j)}_m(\phi) \}\\
&=\sum_{m=1}^{j}\Bigl(\bigl\lfloor\alpha^{(j)}_m(\phi)\bigr\rfloor+m\Bigr),  
\end{align*}
where we used $\alpha^{(j)}_m(\phi)>-m$, see Lemma~\ref{cor:roots}\,(i). This is the first identity in \eqref{eq:landau}. The second one follows by Vieta's formula \eqref{eq:vieta_alpha} from Lemma~\ref{cor:roots}\,(iii). Indeed,
\[
\sum_{m=1}^{j}\Bigl(\bigl\lfloor\alpha^{(j)}_m(\phi)\bigr\rfloor+m\Bigr)
=\sum_{m=1}^{j}\alpha^{(j)}_m(\phi)+\frac{j(j+1)}{2}-\sum_{m=1}^{j}\bigl\{\alpha^{(j)}_m(\phi)\bigr\}
=j\phi-\sum_{m=1}^{j}\bigl\{\alpha^{(j)}_m(\phi)\bigr\}.
\]
Finally, \eqref{eq:twosided} together with the characterization of equality follows from
$0\le\{t\}<1$. 
\end{proof}

\begin{proof}[Proof of Corollary~\ref{cor:flux}]
(i) We have $n\le N(\lambda_n(D_R,B),H_B^{D_R})\le N((2k+1)B,H_B^{D_R})\le k\phi$. Thus, since $B=2\phi/R^2$ and $|D_R|=\pi R^2$,
\[
\lambda_n(D_R,B)>(2k-1)B=\Bigl(1-\frac1{2k}\Bigr)\frac{4\pi}{|D_R|}\,k\phi\ge\Bigl(1-\frac{1}{2k}\Bigr)\frac{4\pi n}{|D_R|}.
\]
(ii) Let $\lambda>0$ and $j=\#\{k\in\N:(2k-1)B<\lambda\}$. Then $\lambda\le(2j+1)B$, and by \eqref{eq:twosided} (for $j\ge1$) and
$N(B,H_B^{D_R})=0$ (for $j=0$),
\[
\#\{n:\lambda_n(D_R,B)<\lambda\}\le N((2j+1)B,H_B^{D_R})\le j\phi .
\]  
\end{proof}

\begin{remark}\label{rem:Bto0}
Fix $n\in\N$. By continuity of $B \mapsto \lambda_n(D_R,B)$, we have $\lambda_n(D_R,B)\to\lambda_n(D_R,0)>0$ as $B\to0$. Hence the integer $k=k(B)$ determined by
$(2k-1)B<\lambda_n(D_R,B)\le(2k+1)B$ tends to infinity as $B\to 0$, and Corollary~\ref{cor:flux}\,(i) implies
$$\lambda_n(D_R,0)\ge \frac{4\pi n}{|D_R|}, \qquad n\in\N .$$
This derivation of P\'olya's inequality for the disk only uses the upper bound in
\eqref{eq:twosided}. In Section~\ref{sec:weyl} we obtain the strict inequality together with an exact formula for the remainder.
\end{remark}

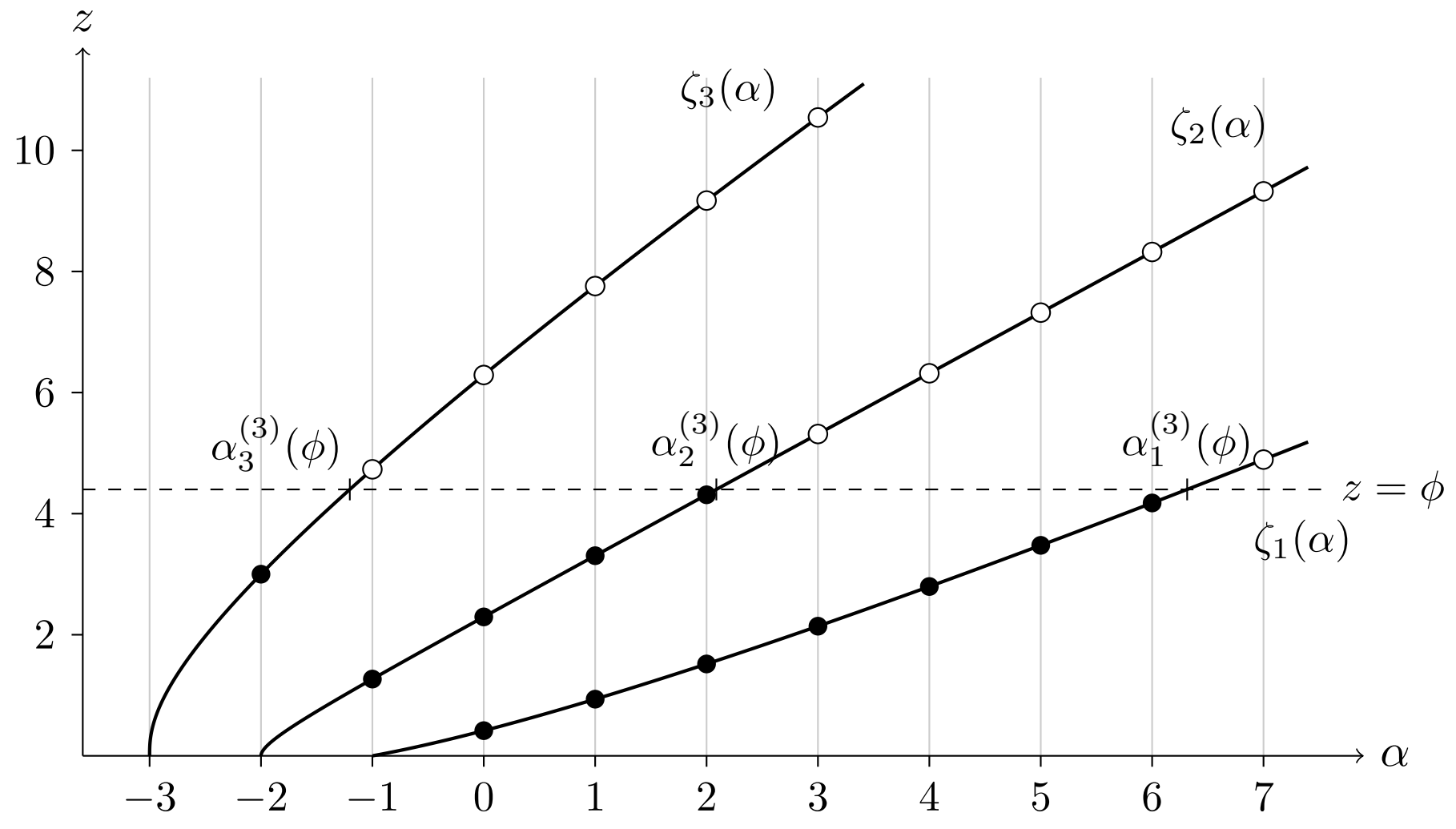
\begin{figure}[t]
\centering
\begin{tikzpicture}[x=0.92cm,y=0.5cm]
\clip (-4.5,-1.3) rectangle (9.0,12.3);
\foreach \a in {-3,...,7} \draw[gray!40] (\a,0) -- (\a,11.2);
\draw[->] (-3.6,0) -- (7.9,0) node[right] {$\alpha$};
\draw[->] (-3.6,0) -- (-3.6,11.7) node[above] {$z$};
\foreach \a in {-3,...,7} \draw (\a,0) -- (\a,-0.18) node[below,font=\small] {$\a$};
\foreach \y in {2,4,6,8,10} \draw (-3.6,\y) -- (-3.7,\y) node[left,font=\small] {$\y$};
\draw[dashed] (-3.6,4.4) -- (7.55,4.4) node[right,font=\small] {$z=\phi$};
\draw[thick] plot coordinates {(-1.0000,0.0000) (-1.0000,0.0000) (-0.9999,0.0000) (-0.9995,0.0002) (-0.9989,0.0004) (-0.9979,0.0007) (-0.9964,0.0012) (-0.9942,0.0019) (-0.9914,0.0029) (-0.9877,0.0041) (-0.9831,0.0057) (-0.9776,0.0075) (-0.9709,0.0098) (-0.9630,0.0125) (-0.9537,0.0157) (-0.9431,0.0193) (-0.9309,0.0235) (-0.9172,0.0283) (-0.9017,0.0338) (-0.8844,0.0400) (-0.8651,0.0469) (-0.8439,0.0546) (-0.8205,0.0632) (-0.7949,0.0727) (-0.7670,0.0832) (-0.7366,0.0948) (-0.7037,0.1075) (-0.6682,0.1215) (-0.6299,0.1367) (-0.5889,0.1534) (-0.5448,0.1715) (-0.4978,0.1911) (-0.4476,0.2124) (-0.3942,0.2354) (-0.3374,0.2603) (-0.2772,0.2871) (-0.2135,0.3159) (-0.1461,0.3469) (-0.0750,0.3802) (0.0000,0.4158) (0.0747,0.4518) (0.1495,0.4884) (0.2242,0.5255) (0.2990,0.5630) (0.3737,0.6011) (0.4485,0.6396) (0.5232,0.6785) (0.5980,0.7178) (0.6727,0.7575) (0.7475,0.7976) (0.8222,0.8381) (0.8970,0.8790) (0.9717,0.9201) (1.0465,0.9617) (1.1212,1.0035) (1.1960,1.0457) (1.2707,1.0882) (1.3455,1.1310) (1.4202,1.1741) (1.4949,1.2175) (1.5697,1.2611) (1.6444,1.3050) (1.7192,1.3492) (1.7939,1.3936) (1.8687,1.4383) (1.9434,1.4832) (2.0182,1.5284) (2.0929,1.5738) (2.1677,1.6194) (2.2424,1.6652) (2.3172,1.7113) (2.3919,1.7576) (2.4667,1.8040) (2.5414,1.8507) (2.6162,1.8976) (2.6909,1.9446) (2.7657,1.9919) (2.8404,2.0393) (2.9152,2.0869) (2.9899,2.1347) (3.0646,2.1827) (3.1394,2.2309) (3.2141,2.2792) (3.2889,2.3276) (3.3636,2.3763) (3.4384,2.4251) (3.5131,2.4740) (3.5879,2.5231) (3.6626,2.5724) (3.7374,2.6218) (3.8121,2.6713) (3.8869,2.7210) (3.9616,2.7709) (4.0364,2.8208) (4.1111,2.8709) (4.1859,2.9212) (4.2606,2.9715) (4.3354,3.0220) (4.4101,3.0726) (4.4848,3.1234) (4.5596,3.1743) (4.6343,3.2253) (4.7091,3.2764) (4.7838,3.3276) (4.8586,3.3789) (4.9333,3.4304) (5.0081,3.4820) (5.0828,3.5337) (5.1576,3.5855) (5.2323,3.6373) (5.3071,3.6894) (5.3818,3.7415) (5.4566,3.7937) (5.5313,3.8460) (5.6061,3.8984) (5.6808,3.9509) (5.7556,4.0035) (5.8303,4.0562) (5.9051,4.1091) (5.9798,4.1619) (6.0545,4.2149) (6.1293,4.2680) (6.2040,4.3212) (6.2788,4.3745) (6.3535,4.4278) (6.4283,4.4812) (6.5030,4.5348) (6.5778,4.5884) (6.6525,4.6421) (6.7273,4.6959) (6.8020,4.7497) (6.8768,4.8036) (6.9515,4.8577) (7.0263,4.9118) (7.1010,4.9659) (7.1758,5.0202) (7.2505,5.0745) (7.3253,5.1289) (7.4000,5.1834)};
\draw[thick] plot coordinates {(-2.0000,0.0000) (-2.0000,0.0024) (-1.9999,0.0068) (-1.9995,0.0126) (-1.9989,0.0196) (-1.9979,0.0277) (-1.9964,0.0369) (-1.9942,0.0472) (-1.9914,0.0585) (-1.9877,0.0710) (-1.9831,0.0845) (-1.9776,0.0993) (-1.9709,0.1152) (-1.9630,0.1324) (-1.9537,0.1509) (-1.9431,0.1708) (-1.9309,0.1921) (-1.9172,0.2149) (-1.9017,0.2392) (-1.8844,0.2652) (-1.8651,0.2929) (-1.8439,0.3223) (-1.8205,0.3537) (-1.7949,0.3870) (-1.7670,0.4223) (-1.7366,0.4598) (-1.7037,0.4996) (-1.6682,0.5417) (-1.6299,0.5862) (-1.5889,0.6333) (-1.5448,0.6831) (-1.4978,0.7356) (-1.4476,0.7909) (-1.3942,0.8493) (-1.3374,0.9107) (-1.2772,0.9753) (-1.2135,1.0432) (-1.1461,1.1146) (-1.0750,1.1894) (-1.0000,1.2679) (-0.9152,1.3564) (-0.8303,1.4444) (-0.7455,1.5321) (-0.6606,1.6196) (-0.5758,1.7068) (-0.4909,1.7938) (-0.4061,1.8806) (-0.3212,1.9673) (-0.2364,2.0538) (-0.1515,2.1402) (-0.0667,2.2265) (0.0182,2.3127) (0.1030,2.3989) (0.1879,2.4849) (0.2727,2.5708) (0.3576,2.6567) (0.4424,2.7426) (0.5273,2.8283) (0.6121,2.9141) (0.6970,2.9997) (0.7818,3.0854) (0.8667,3.1710) (0.9515,3.2565) (1.0364,3.3421) (1.1212,3.4275) (1.2061,3.5130) (1.2909,3.5984) (1.3758,3.6839) (1.4606,3.7692) (1.5455,3.8546) (1.6303,3.9400) (1.7152,4.0253) (1.8000,4.1106) (1.8848,4.1959) (1.9697,4.2811) (2.0545,4.3664) (2.1394,4.4516) (2.2242,4.5369) (2.3091,4.6221) (2.3939,4.7073) (2.4788,4.7925) (2.5636,4.8776) (2.6485,4.9628) (2.7333,5.0479) (2.8182,5.1331) (2.9030,5.2182) (2.9879,5.3034) (3.0727,5.3885) (3.1576,5.4736) (3.2424,5.5587) (3.3273,5.6438) (3.4121,5.7289) (3.4970,5.8140) (3.5818,5.8990) (3.6667,5.9841) (3.7515,6.0692) (3.8364,6.1542) (3.9212,6.2393) (4.0061,6.3243) (4.0909,6.4094) (4.1758,6.4944) (4.2606,6.5794) (4.3455,6.6645) (4.4303,6.7495) (4.5152,6.8345) (4.6000,6.9195) (4.6848,7.0045) (4.7697,7.0895) (4.8545,7.1745) (4.9394,7.2595) (5.0242,7.3445) (5.1091,7.4295) (5.1939,7.5145) (5.2788,7.5995) (5.3636,7.6845) (5.4485,7.7695) (5.5333,7.8544) (5.6182,7.9394) (5.7030,8.0244) (5.7879,8.1094) (5.8727,8.1943) (5.9576,8.2793) (6.0424,8.3643) (6.1273,8.4492) (6.2121,8.5342) (6.2970,8.6191) (6.3818,8.7041) (6.4667,8.7890) (6.5515,8.8740) (6.6364,8.9589) (6.7212,9.0439) (6.8061,9.1288) (6.8909,9.2138) (6.9758,9.2987) (7.0606,9.3837) (7.1455,9.4686) (7.2303,9.5535) (7.3152,9.6385) (7.4000,9.7234)};
\draw[thick] plot coordinates {(-3.0000,0.0000) (-3.0000,0.0328) (-2.9999,0.0668) (-2.9995,0.1021) (-2.9989,0.1386) (-2.9979,0.1767) (-2.9964,0.2162) (-2.9942,0.2574) (-2.9914,0.3003) (-2.9877,0.3451) (-2.9831,0.3917) (-2.9776,0.4404) (-2.9709,0.4912) (-2.9630,0.5442) (-2.9537,0.5994) (-2.9431,0.6571) (-2.9309,0.7173) (-2.9172,0.7800) (-2.9017,0.8454) (-2.8844,0.9136) (-2.8651,0.9846) (-2.8439,1.0585) (-2.8205,1.1356) (-2.7949,1.2157) (-2.7670,1.2990) (-2.7366,1.3857) (-2.7037,1.4758) (-2.6682,1.5694) (-2.6299,1.6666) (-2.5889,1.7675) (-2.5448,1.8722) (-2.4978,1.9807) (-2.4476,2.0933) (-2.3942,2.2098) (-2.3374,2.3306) (-2.2772,2.4556) (-2.2135,2.5849) (-2.1461,2.7187) (-2.0750,2.8570) (-2.0000,3.0000) (-1.9051,3.1774) (-1.8101,3.3511) (-1.7152,3.5217) (-1.6202,3.6894) (-1.5253,3.8546) (-1.4303,4.0175) (-1.3354,4.1783) (-1.2404,4.3372) (-1.1455,4.4944) (-1.0505,4.6499) (-0.9556,4.8040) (-0.8606,4.9567) (-0.7657,5.1080) (-0.6707,5.2582) (-0.5758,5.4072) (-0.4808,5.5552) (-0.3859,5.7021) (-0.2909,5.8481) (-0.1960,5.9932) (-0.1010,6.1374) (-0.0061,6.2808) (0.0889,6.4235) (0.1838,6.5653) (0.2788,6.7065) (0.3737,6.8470) (0.4687,6.9869) (0.5636,7.1262) (0.6586,7.2648) (0.7535,7.4029) (0.8485,7.5404) (0.9434,7.6774) (1.0384,7.8139) (1.1333,7.9499) (1.2283,8.0854) (1.3232,8.2205) (1.4182,8.3551) (1.5131,8.4892) (1.6081,8.6230) (1.7030,8.7564) (1.7980,8.8894) (1.8929,9.0219) (1.9879,9.1542) (2.0828,9.2860) (2.1778,9.4176) (2.2727,9.5487) (2.3677,9.6796) (2.4626,9.8101) (2.5576,9.9404) (2.6525,10.0703) (2.7475,10.1999) (2.8424,10.3292) (2.9374,10.4583) (3.0323,10.5871) (3.1273,10.7156) (3.2222,10.8438) (3.3172,10.9718) (3.4121,11.0996)};
\fill (0,0.4158) circle (2.2pt);
\fill (1,0.9358) circle (2.2pt);
\fill (2,1.5174) circle (2.2pt);
\fill (3,2.1412) circle (2.2pt);
\fill (4,2.7965) circle (2.2pt);
\fill (5,3.4764) circle (2.2pt);
\fill (6,4.1763) circle (2.2pt);
\draw[fill=white] (7,4.8928) circle (2.2pt);
\fill (-1,1.2679) circle (2.2pt);
\fill (0,2.2943) circle (2.2pt);
\fill (1,3.3054) circle (2.2pt);
\fill (2,4.3116) circle (2.2pt);
\draw[fill=white] (3,5.3155) circle (2.2pt);
\draw[fill=white] (4,6.3182) circle (2.2pt);
\draw[fill=white] (5,7.3202) circle (2.2pt);
\draw[fill=white] (6,8.3218) circle (2.2pt);
\draw[fill=white] (7,9.3230) circle (2.2pt);
\fill (-2,3.0000) circle (2.2pt);
\draw[fill=white] (-1,4.7321) circle (2.2pt);
\draw[fill=white] (0,6.2899) circle (2.2pt);
\draw[fill=white] (1,7.7588) circle (2.2pt);
\draw[fill=white] (2,9.1710) circle (2.2pt);
\draw[fill=white] (3,10.5433) circle (2.2pt);
\draw (6.3146,4.22) -- (6.3146,4.58) node[above=-1pt,font=\small] {$\alpha^{(3)}_1(\phi)$};
\draw (2.0880,4.22) -- (2.0880,4.58) node[above=-1pt,font=\small] {$\alpha^{(3)}_2(\phi)$};
\draw (-1.2026,4.22) -- (-1.2026,4.58) node[above left=-2pt,font=\small] {$\alpha^{(3)}_3(\phi)$};
\node[font=\small] at (7.35,3.55) {$\zeta_1(\alpha)$};
\node[font=\small] at (6.6,10.4) {$\zeta_2(\alpha)$};
\node[font=\small] at (2.2,11.0) {$\zeta_3(\alpha)$};
\end{tikzpicture}
\caption{An illustration of Theorem~\ref{thm:landau}. For $j=3$, the zero set $Z$ consists of the zero curves $\zeta_1,\zeta_2,\zeta_3$ of Proposition~\ref{prop:curves}. The curve $\zeta_m$ emanates from the point $(-m,0)$ and is strictly increasing with $\alpha$. }
\label{fig:curves2}
\end{figure}

\begin{remark}
Figure~\ref{fig:curves2} illustrates the two identities of Theorem~\ref{thm:landau} for $j=3$. The dots are the points of the curves $\zeta_m(\alpha)$ with integer $\alpha$. A point $(\alpha,\zeta_m(\alpha))$ with $\alpha=-l\in\Z$ corresponds to an eigenvalue branch with angular momentum $l$ which crosses the Landau level $7B$ at flux $\phi=\zeta_m(\alpha)$. In the picture, there are $12$ points (filled dots) on or below the dashed line at $\phi=4.4$, so by \eqref{eq:Nall}, $N(7B,H_B^{D_R})=12$. The marks on the dashed line are the zeros $\alpha^{(3)}_m(\phi)\approx6.31,\,2.09,\,-1.20$ from Lemma~\ref{cor:roots}. Then the number of filled dots on each curve $\zeta_m(\alpha)$ is precisely $(\lfloor\alpha^{(3)}_m(\phi)\rfloor+m)$, and indeed
\begin{align*}
    \sum_{m=1}^{3}\Bigl(\bigl\lfloor\alpha^{(3)}_m(\phi)\bigr\rfloor+m\Bigr) &= (6+1)+(2+2)+(-2+3)=12 
\end{align*}
and
\begin{align*}
    3\phi-\sum_{m=1}^{3}\bigl\{\alpha^{(3)}_m(\phi)\bigr\} &= 3 \cdot 4.4-(0.31+0.09+0.80) = 12.
\end{align*}
\end{remark}

\section{Zeros of Bessel functions as functions of the order}\label{sec:bessel}

This section runs largely parallel to Section~\ref{sec:laguerre}. The function $E$ defined below takes the role of $P_j$, and the
order $\nu$ takes the role of the parameter $\alpha$. Lemma~\ref{lem:E}, Proposition~\ref{prop:besselcurves} and
Lemma~\ref{thm:bessel} correspond to Lemma~\ref{lem:Pj}, Proposition~\ref{prop:curves} and Lemma~\ref{cor:roots},
respectively. The counterpart of Vieta's formula \eqref{eq:vieta_alpha} of
Lemma~\ref{cor:roots}\,(iii) is the sum rule \eqref{eq:sumrule} of Lemma~\ref{thm:bessel}\,(iii).

For $(\nu,x)\in\C^2$ let
\begin{equation}\label{eq:E}
E(\nu,x):=\sum_{k=0}^{\infty}\frac{(-\tfrac{x^2}{4})^k}{k!\,\Gamma(\nu+k+1)},
\end{equation}
where $1/\Gamma$ is understood as an entire function. By the power series of the Bessel function \cite[\S10.2(ii)]{DLMF},
\begin{equation}\label{eq:EJ}
E(\nu,x)=\Bigl(\frac x2\Bigr)^{-\nu}J_\nu(x)\qquad\text{for }x>0,\ \nu\in\C .
\end{equation}
Hence, for $x>0$, the zeros of $\nu\mapsto J_\nu(x)$ are the zeros of $E(\cdot,x)$, and for $\nu\in\R$ the positive zeros of $J_\nu$
are those of $E(\nu,\cdot)$. The link to Section~\ref{sec:laguerre} is provided by the functions
\begin{equation}\label{eq:Ej}
E_j(\nu,x):=\frac{P_j\bigl(\nu,\tfrac{x^2}{4j}\bigr)}{\Gamma(\nu+j+1)}=\sum_{i=0}^{j}\binom ji\,j^{-i}\,\frac{(-\tfrac{x^2}{4})^i}{\Gamma(\nu+i+1)},
\qquad j\in\N,
\end{equation}
where $P_j$ is extended to $\C^2$ by \eqref{eq:Pj}. The second expression follows from \eqref{eq:Pj} and the fact that
$(\nu+i+1)_{j-i}\cdot 1/\Gamma(\nu+j+1)= 1/\Gamma(\nu+i+1)$.

\begin{lemma}\label{lem:E} \,
\begin{enumerate}[label=\textup{(\roman*)}]
\item $E$ is an entire function on $\C^2$ which is real for real arguments, and $E_j\to E$ locally uniformly on $\C^2$ as $j\to\infty$.
\item For every $\nu\in\R$ the function $y=E(\nu,\cdot)$ satisfies $xy''+(2\nu+1)y'+xy=0$ on $(0,\infty)$.
Consequently, all zeros of $E(\nu,\cdot)$ in $(0,\infty)$ are simple.
\item $E(\nu,0)=1/\Gamma(\nu+1)$, which vanishes if and only if $\nu\in\{-1,-2,...\}$.
\item If $\nu\in\R$, then $E(\nu,\cdot)$ has infinitely many zeros in $(0,\infty)$. If $x>0$ and $\nu\in\R$ with
$\nu+1>\frac{x^2}{4\log2}$, then $E(\nu,x)>0$ and $E_j(\nu,x)>0$ for all $j\in\N$.
\item For $k\in\N$ and $x\in\C$,
\begin{equation}\label{eq:negorder}
E(-k,x)=\Bigl(-\frac{x^2}{4}\Bigr)^{k}E(k,x).
\end{equation}
In particular, for every $l\in\Z$ the zeros of $E(l,\cdot)$ in $(0,\infty)$ are the positive zeros of $J_{|l|}$.
\end{enumerate}
\end{lemma}

\begin{proof}
(i) Let $r>0$ and $k_0:=\lceil r\rceil+1$. For $|\nu|\le r$ and $k>k_0$ we have
$1/\Gamma(\nu+k+1)= 1/\Gamma(\nu+k_0+1) \cdot (\prod_{i=k_0+1}^{k}(\nu+i))^{-1}$ and $|\nu+i|\ge1$ for $i>k_0$. Therefore, $| 1/\Gamma(\nu+k+1)| \leq |1/ \Gamma(\nu+k_0+1)|$ for $k>k_0$ and hence the series \eqref{eq:E} is dominated by $C_r\sum_k(|x|^2/4)^k/k!$ for $|\nu|\le r$ where $C_r:=\sup_{k\in\N_0}\sup_{|\nu|\le r}| 1/\Gamma(\nu+k+1)|<\infty$. Thus, the series converges absolutely and uniformly for $|\nu|\le r$, $|x|\le r$ and since $r>0$ was arbitrary, $E$ is entire in $\nu$ and $x$. Since $0\le\binom jij^{-i}\le\frac1{i!}$ and
$\binom jij^{-i}\to\frac1{i!}$ as $j\to\infty$ for fixed $i$, the same bound and dominated convergence give $E_j\to E$ uniformly for
$|\nu|\le r$, $|x|\le r$.

(ii) By \eqref{eq:EJ}, the differential equation is equivalent to Bessel's equation $x^2w''+xw'+(x^2-\nu^2)w=0$ for $w=J_\nu$, see
\cite[(10.2.1)]{DLMF}. It is regular on $(0,\infty)$, and simplicity of the zeros in $(0,\infty)$ follows as in
Lemma~\ref{lem:Pj}\,(ii).

(iii) is immediate from \eqref{eq:E}.

(iv) For real $\nu$ the function $J_\nu$ has infinitely many positive zeros \cite[\S10.21(i)]{DLMF}. For the second statement, note
that by \eqref{eq:E} and \eqref{eq:Ej},
\begin{align*}
  \Gamma(\nu+1)E(\nu,x)&=\sum_{k=0}^\infty\frac{(-\tfrac{x^2}{4})^k}{k!\,(\nu+1)_k}, \\ 
  \Gamma(\nu+1)E_j(\nu,x)&=\sum_{i=0}^{j}\binom jij^{-i}\frac{(-\tfrac{x^2}{4})^i}{(\nu+1)_i}.
\end{align*}
As $(\nu+1)_k\ge(\nu+1)^k$, both sums differ from $1$ by at most 
$$\sum_{k=1}^\infty\frac{1}{k!}\Bigl(\frac{x^2}{4(\nu+1)}\Bigr)^k=e^{x^2/(4(\nu+1))}-1<1$$ 
for $\nu+1>\frac{x^2}{4\log2}$. Thus, $E(\nu,x)>0$ and $E_j(\nu,x)>0$ for all $j\in\N$ and $\nu+1>\frac{x^2}{4\log2}$.

(v) For $\nu=-k$ the terms in \eqref{eq:E} with index less than $k$ vanish, and shifting the index of summation by $k$ gives
\eqref{eq:negorder}. The last assertion follows from \eqref{eq:negorder} and \eqref{eq:EJ}.
\end{proof}

The identity \eqref{eq:negorder} is the classical relation $J_{-k}=(-1)^kJ_k$. For fixed $x>0$, the function $E(\cdot,x)$ is not a polynomial in $\nu$, so there is no direct analogue of Lemma~\ref{lem:Pj}\,(i). The following lemma serves as a substitute: instead of the degree and the second coefficient of a polynomial, it provides the number of zeros of $E(\cdot,x)$ in large disks and their sum.

\begin{lemma}\label{lem:circle}
Fix $x>0$ and let
\[
G(\nu):=\Gamma(\nu+1) E(\nu,x)=\sum_{k=0}^{\infty}\frac{(-\tfrac{x^2}{4})^k}{k!\,(\nu+1)_k},\qquad\nu\in\C\setminus\{-1,-2,...\}.
\]
There are constants $N_0(x)\ge3$ and $C(x)$, which can be chosen non-decreasing in $x$, such that for all integers
$N\ge N_0(x)$ the following holds: the function $E(\cdot,x)$ has no zeros on the circle $|\nu|=N+\frac12$ and it has exactly $N$ zeros in the disk $|\nu|<N+\frac12$, counted with multiplicity. Denoting these zeros by $\{\nu_{m}(x)\}_{m=1}^N$, the sum of these zeros satisfies
\begin{equation}\label{eq:sigmaN}
\bigg| \sum_{m=1}^N (\nu_{m}(x) +m)-\frac{x^2}{4}\bigg|\le\frac{C(x)}{N} .
\end{equation}
\end{lemma}

\begin{proof}
We first show that for every integer $N\ge3$ and all $\nu\in\C$ with $|\nu|=N+\frac12$,
\begin{equation}\label{eq:Gcircle}
\Bigl|G(\nu)-1+\frac{x^2}{4(\nu+1)}\Bigr|\le\frac{4e^{x^2/4}}{(N-\frac12)(N-\frac32)} .
\end{equation}
Let $|\nu|=N+\frac12$. By definition,
\[
\Bigl|G(\nu)-1+\frac{x^2}{4(\nu+1)}\Bigr|\le\sum_{k=2}^\infty\frac{(x^2/4)^k}{k!\, |(\nu+1)_k| }.
\]
For $i\in\N$ we have $|\nu+i|\ge |N+\frac12-i|$ and therefore
$|(\nu+1)_k|\ge \prod_{i=1}^{k}|N+\frac12-i|$. Among the factors
$|N+\frac12-i|$ with $i\ge3$ exactly two are smaller than $1$, namely those with $i=N$ and $i=N+1$, which are equal to $\frac12$. Hence, for $k\ge2$,
\begin{align*}
 |(\nu+1)_k|\ge \frac{1}{4} \prod_{i=1}^{2}|N+\tfrac12-i| = \frac{1}{4} (N-\tfrac12)(N-\tfrac32)
\end{align*}
and
\begin{align*}
    \Bigl|G(\nu)-1+\frac{x^2}{4(\nu+1)}\Bigr|\le\frac{4}{(N-\frac12)(N-\frac32)}\sum_{k =2}^\infty\frac{(x^2/4)^k}{k!}\le\frac{4e^{x^2/4}}{(N-\frac12)(N-\frac32)}.
\end{align*}

Let now
\begin{align*}
    C^{(1)}_N=\frac{x^2}{4(N-\frac12)}, \qquad C^{(2)}_N=\frac{4e^{x^2/4}}{(N-\frac12)(N-\frac32)}
\end{align*}
and choose $N_0(x)\ge3$ such that $C^{(1)}_N+C^{(2)}_N\le\frac12$ for
$N\ge N_0(x)$. Since $C^{(1)}_N+C^{(2)}_N$ is increasing in $x$, $N_0(x)$ can be chosen non-decreasing in $x$. Let $N\ge N_0(x)$ and let $\mathbb{T}_{N+1/2}$ denote the positively oriented circle $|\nu|=N+\frac12$. By \eqref{eq:Gcircle} and $|\nu+1|\ge N-\frac12$, we have $|G(\nu)-1|\le C^{(1)}_N+C^{(2)}_N \le\frac12$
on $\mathbb{T}_{N+1/2}$. In particular $G$, and hence $E(\cdot,x)$, has no zeros on $\mathbb{T}_{N+1/2}$, and $\log G$, defined with the principal branch of the
logarithm, is holomorphic in a neighbourhood of $\mathbb{T}_{N+1/2}$ with $(\log G)'=G'/G$. Using $|\log(1+w)-w|\le|w|^2$ for $|w|\le\frac12$ with $w=G(\nu)-1$ and
\eqref{eq:Gcircle} again, we obtain for $\nu\in\C$ with $|\nu|=N+\frac12$
\begin{align}\label{eq:logG}
\begin{split}
\Bigl|\log G(\nu)+\frac{x^2}{4(\nu+1)}\Bigr|& \leq \Bigl|\log G(\nu)-(G(\nu)-1)\Bigr| + \Bigl|G(\nu)-1+\frac{x^2}{4(\nu+1)}\Bigr| \\
&\le (C^{(1)}_N+C^{(2)}_N)^2 + C^{(2)}_N.
\end{split}
\end{align}
The function $G$ is meromorphic in $\C$. Since $1/\Gamma(\nu+1)$ is entire with simple zeros exactly at the negative integers, the order of $G$ at a point $\nu_0$ is the order of $E(\cdot,x)$ at $\nu_0$, diminished by $1$ if $\nu_0\in\{-1,-2,...\}$. Let
$Z_N$ denote the number of zeros of $E(\cdot,x)$ inside $\mathbb{T}_{N+1/2}$, counted with multiplicity. By the argument principle
\cite[Ch.~V, \S3]{Conway}, applied to the functions $1$ and $\nu$,
\[
Z_N-N=\frac{1}{2\pi i}\oint_{\mathbb{T}_{N+1/2}}\frac{G'(\nu)}{G(\nu)}\,\dd\nu,
\qquad
\sum_{m=1}^N \nu_{m}(x)-\sum_{m=1}^{N}(-m)=\frac{1}{2\pi i}\oint_{\mathbb{T}_{N+1/2}}\nu\,\frac{G'(\nu)}{G(\nu)}\,\dd\nu .
\]
The first integral is zero because $G'/G$ has the antiderivative $\log G$ along $\mathbb{T}_{N+1/2}$, so $Z_N=N$. The second integral we integrate
by parts and use \eqref{eq:logG} and $\frac{1}{2\pi i}\oint_{\mathbb{T}_{N+1/2}}\frac{\dd\nu}{\nu+1}=1$ to obtain
\[
\sum_{m=1}^N (\nu_{m}(x)+m)=-\frac{1}{2\pi i}\oint_{\mathbb{T}_{N+1/2}}\log G(\nu)\,\dd\nu=\frac{x^2}{4}+R_N
\]
with
\[
|R_N|\le (N+\tfrac12)((C^{(1)}_N+C^{(2)}_N)^2 + C^{(2)}_N).
\]
Since $C^{(1)}_N=O(N^{-1})$ and $C^{(2)}_N=O(N^{-2})$ as $N\to \infty$, this shows \eqref{eq:sigmaN} and since the right-hand side is increasing in $x$, $C(x)$ can be chosen non-decreasing in $x$.
\end{proof}

The next proposition is the counterpart of Proposition~\ref{prop:curves} for zeros of Bessel functions. The only difference in the statement is that there are
now infinitely many curves.

\begin{proposition}\label{prop:besselcurves}
Let
$$Z:=\{(\nu,x)\in\R\times(0,\infty):E(\nu,x)=0\}.$$
There exist continuous, strictly increasing bijections
\begin{align*}
    x_m\colon(-m,\infty)&\to(0,\infty),\qquad m\in\N, \\
\nu &\mapsto x_m(\nu),
\end{align*}
with the following properties:
\begin{enumerate}[label=\textup{(\roman*)}]
\item $Z=\bigcup_{m=1}^{\infty}\bigl\{(\nu,x_m(\nu)):\nu>-m\bigr\}$.
\item If $m,m'\in\N$ with $m<m'$, then $x_m(\nu)<x_{m'}(\nu)$ for $\nu>-m$.
\end{enumerate}
In other words, the zero set $Z$ consists precisely of countably many non-intersecting, strictly increasing graphs.
\end{proposition}

By (i) and (ii), for $\nu>-1$ the number $x_m(\nu)$ is the $m$-th positive zero of $J_\nu$, which is usually denoted by $j_{\nu,m}$.
Thus $x_m$ is the continuation of $\nu\mapsto j_{\nu,m}$ to the interval $(-m,\infty)$, and as in Proposition~\ref{prop:curves} the
curve $x_m$ emanates from the point $(-m,0)$. That $j_{\nu,m}$ can be continued to an increasing function of $\nu$ on
$(-m,\infty)$ which tends to $0$ as $\nu\searrow-m$ is known, see Elbert \cite{Elbert77} and the survey \cite[Sections~1.1 and~1.3]{Elbert01}.
The proof below uses the same uniqueness principle as the proof of
Proposition~\ref{prop:curves}, but the curves are not constructed by continuation in the parameter. Instead, we first construct their
inverse functions $\nu_m$ as limits of the functions $\alpha^{(j)}_m$ from Lemma~\ref{cor:roots}.

\begin{proof}
By Lemma~\ref{lem:E}\,(ii), $\partial_xE\neq0$ on $Z$. By the implicit function theorem,
around every point of $Z$, $Z$ is locally the graph of a continuous function of $\nu$. As in the proof of
Proposition~\ref{prop:curves}, this implies the following uniqueness principle: if $I$ is an interval and
$\xi,\tilde{\xi}\colon I\to(0,\infty)$ are continuous functions whose graphs are contained in $Z$ and which coincide at one point, then
$\xi=\tilde{\xi}$ on $I$.

\emph{Step 1: Construct $\nu_m$.} Fix $x>0$. By \eqref{eq:Ej} and Lemma~\ref{cor:roots}\,(i), the zeros of the entire function
$E_j(\cdot,x)$ are the real numbers
\begin{equation}\label{eq:Ejzeros}
\alpha^{(j)}_1\bigl(\tfrac{x^2}{4j}\bigr)>\alpha^{(j)}_2\bigl(\tfrac{x^2}{4j}\bigr)>...>\alpha^{(j)}_j\bigl(\tfrac{x^2}{4j}\bigr),
\end{equation}
together with the zeros $-j-1,-j-2,...$ of $1/\Gamma(\nu+j+1)$. By Lemma~\ref{lem:E}\,(i) and (iv),
$E_j(\cdot,x)\to E(\cdot,x)\not\equiv0$ locally uniformly. Hence, by Hurwitz's theorem \cite[Ch.~VII, \S2]{Conway}, every zero of
$E(\cdot,x)$ is a limit of zeros of $E_j(\cdot,x)$ and therefore real. Now let $N_0(x)$ be
the constant from Lemma~\ref{lem:circle} and let $N\ge N_0(x)$ be an integer with $N+\frac12>\frac{x^2}{4\log2}$.
By Lemma~\ref{lem:circle} and Lemma~\ref{lem:E}\,(iv), the function $E(\cdot,x)$ has exactly $N$ zeros in the interval
$(-N-\frac12,N+\frac12)$, counted with multiplicity, and no zeros in $[N+\frac12,\infty)$. Since $N$ can be chosen arbitrarily large,
$E(\cdot,x)$ has infinitely many zeros, and we can enumerate them in non-increasing order and repeated according to multiplicity as
$\nu_1(x)\ge\nu_2(x)\ge...$. Then, for every $N$ as above, the zeros of $E(\cdot,x)$ in the disk $|\nu|<N+\frac12$, counted
with multiplicity, are the numbers
\begin{equation}\label{eq:firstN}
\nu_1(x)\ge\nu_2(x)\ge...\ge\nu_N(x).
\end{equation}

\emph{Step 2: $\nu_m$ as a limit of Laguerre zeros.} We claim that
\begin{equation}\label{eq:limit}
\nu_m(x)=\lim_{j\to\infty}\alpha^{(j)}_m\Bigl(\frac{x^2}{4j}\Bigr)\qquad\text{for all }x>0,\ m\in\N .
\end{equation}
Fix $x>0$, let $N_0(x)$ be the constant from Lemma~\ref{lem:circle} and $m\in\N$. Choose $N\ge \max\{m,N_0(x)\}$ with $N+\frac12>\frac{x^2}{4\log2}$. For $j\ge N$, the first $N$ numbers in \eqref{eq:Ejzeros} are
larger than $-N$ by Lemma~\ref{cor:roots}\,(i) and smaller than $N+\frac12$ by Lemma~\ref{lem:E}\,(iv), so they lie inside the disk $|\nu|<N+\frac12$. Since $E(\cdot,x)$ has no zeros on the circle $|\nu|=N+\frac12$ and exactly $N$ zeros inside by Lemma~\ref{lem:circle}, Hurwitz's theorem
shows that for all sufficiently large $j$ the function $E_j(\cdot,x)$ has exactly $N$ zeros in this disk, which are then the first
$N$ numbers in \eqref{eq:Ejzeros}, and that these zeros converge to the zeros of $E(\cdot,x)$ in the disk, counted with multiplicity.
As both families are listed in non-increasing order, \eqref{eq:limit} follows from \eqref{eq:firstN}.

\emph{Step 3: Monotonicity, continuity and limit of $\nu_m$ as $x\searrow0$.}
Since each $\alpha^{(j)}_m$ is increasing with values in $(-m,\infty)$, \eqref{eq:limit} shows that
\begin{equation}\label{eq:numweak}
\nu_m\text{ is non-decreasing on }(0,\infty)\qquad\text{and}\qquad\nu_m(x)\ge-m\text{ for all }x>0 .
\end{equation}
The same argument applies to $E(\cdot,x')$ as $x'\to x$ instead of $E_j(\cdot,x)$ as $j\to\infty$, because $E(\cdot,x')\to E(\cdot,x)$
locally uniformly and \eqref{eq:firstN} holds with the same $N$ for all $x'$ in a neighbourhood of $x$: if $N$ satisfies the
conditions of Step~1 for some $x''>x$, then it satisfies them for all $x'\in(0,x'']$, as $N_0(x)$ and $\frac{x^2}{4\log2}$ are
non-decreasing in $x$. This shows that each $\nu_m$ is continuous on $(0,\infty)$. Likewise, $E(\cdot,x')\to E(\cdot,0)=1/\Gamma(\cdot+1)$ locally
uniformly as $x'\searrow0$, and by Lemma~\ref{lem:E}\,(iii) the zeros of the limit are the simple zeros $-1,-2,...$. Hence, again
by the same argument,
\begin{equation}\label{eq:nulimit0}
\nu_m(x)\to-m,\qquad x\searrow0 .
\end{equation}

\emph{Step 4: Strict monotonicity and range of $\nu_m$.} We claim that each $\nu_m$ is strictly
increasing. Suppose that $\nu_m(x_1)=\nu_m(x_2)=\nu^*$ for some $0<x_1<x_2$. By
\eqref{eq:numweak}, $\nu_m=\nu^*$ on $[x_1,x_2]$, so $E(\nu^*,\cdot)$ vanishes on $[x_1,x_2]$. This contradicts
Lemma~\ref{lem:E}\,(ii), by which the zeros of $E(\nu^*,\cdot)$ in $(0,\infty)$ are simple. Hence each $\nu_m$ is strictly
increasing.

Next we show that
\begin{equation}\label{eq:nulimitinf}
\nu_m(x)\to\infty,\qquad x\to\infty .
\end{equation}
Suppose that $\nu_m\le\nu^*<\infty$ on $(0,\infty)$, and fix a real number $\nu_0>\nu^*$. By Lemma~\ref{lem:E}\,(iv),
$E(\nu_0,\cdot)$ has infinitely many positive zeros $y_1<y_2<...$. For each $i$ the number $\nu_0$ is a zero of $E(\cdot,y_i)$, so
$\nu_0=\nu_{m_i}(y_i)$ for some $m_i\in\N$, and $m_i<m$ because $\nu_{m_i}(y_i)=\nu_0>\nu_m(y_i)$. Hence $m_i=m_{i'}$ for some
$i\ne i'$, and the function $\nu_{m_i}$ takes the value $\nu_0$ at two different points. This contradicts the strict monotonicity.
By \eqref{eq:nulimit0}, \eqref{eq:nulimitinf} and continuity, each $\nu_m$ is a continuous, strictly increasing bijection from
$(0,\infty)$ onto $(-m,\infty)$.

\emph{Step 5: Construct $x_m$.} Let $x_m:=\nu_m^{-1}\colon(-m,\infty)\to(0,\infty)$. By Step~4, $x_m$ is a continuous, strictly
increasing bijection. Let $\nu\in\R$ and $x>0$. By Step~1, $E(\nu,x)=0$ holds if and only if $\nu=\nu_m(x)$ for some $m\in\N$, that
is, if and only if $\nu>-m$ and $x=x_m(\nu)$ for some $m\in\N$. This is (i).

\emph{Step 6: Ordering of the curves.} Let $m<m'$ and $\nu>-m$. Since $\nu_{m'}\le\nu_m$, we have
$\nu_{m'}(x_m(\nu))\le\nu_m(x_m(\nu))=\nu=\nu_{m'}(x_{m'}(\nu))$, and hence $x_m(\nu)\le x_{m'}(\nu)$ because $\nu_{m'}$ is strictly
increasing. Suppose for the sake of contradiction that we have equality, i.e.~$x_m(\nu_0)=x_{m'}(\nu_0)$ for some $\nu_0>-m$. By (i), the graphs of $x_m$ and $x_{m'}$ are contained in
$Z$, so the uniqueness principle gives $x_m=x_{m'}$ on $(-m,\infty)$. Letting $\nu\searrow-m$, we obtain $x_m(\nu)\to0$, whereas
$x_{m'}(\nu)\to x_{m'}(-m)>0$ since $-m$ lies in the domain $(-m',\infty)$ of $x_{m'}$. This is a contradiction, which proves (ii).
\end{proof}

The following lemma is the counterpart of Lemma~\ref{cor:roots}.

\begin{lemma}\label{thm:bessel}
Let $x>0$, and let $\nu_m\colon(0,\infty)\to(-m,\infty)$, $m\in\N$, be the functions from the proof of Proposition~\ref{prop:besselcurves}.
\begin{enumerate}[label=\textup{(\roman*)}]
\item The zeros of the function $\nu\mapsto J_\nu(x)$, $\nu\in\C$, are exactly the numbers $\nu_m(x)$, $m\in\N$. They are real and
simple, and
\[
\nu_1(x)>\nu_2(x)>\nu_3(x)>...,\qquad\nu_m(x)>-m .
\]
Each $\nu_m$ is continuous and strictly increasing.
\item For every $\nu\in\R$,
\[
\#\{t\in(0,x]:J_\nu(t)=0\}=\#\{m\in\N:-m<\nu\le\nu_m(x)\}.
\]
\item We have
\begin{equation}\label{eq:sumrule}
\sum_{m=1}^{\infty} (\nu_m(x)+m )=\frac{x^2}{4} .
\end{equation}
\end{enumerate}
\end{lemma}

We are not aware of a reference for the sum rule \eqref{eq:sumrule}. That the zeros of $\nu\mapsto J_\nu(x)$ are real, simple and asymptotically close to the negative integers goes back to Coulomb
\cite{Coulomb}, see also \cite[\S10.21(xiv)]{DLMF} and the references given there, as well as \cite{CondeKalla,FlajoletSchott}. Our proof of (i) is independent of these works.

\begin{proof}
By \eqref{eq:EJ}, the zeros of $\nu\mapsto J_\nu(x)$ are the zeros of $E(\cdot,x)$, and for $\nu\in\R$ the zeros of
$J_\nu$ in $(0,\infty)$ are those of $E(\nu,\cdot)$.

(i) By Proposition~\ref{prop:besselcurves}\,(i), for $\nu\in\R$ we have $E(\nu,x)=0$ if and only if $x=x_m(\nu)$ for some $m\in\N$
with $\nu>-m$, or equivalently, $\nu=\nu_m(x)$. If $m<m'$, then $\nu_m(x)>-m>-m'$ lies in the domain of $x_{m'}$, and
$x_{m'}(\nu_m(x))>x_m(\nu_m(x))=x=x_{m'}(\nu_{m'}(x))$ by Proposition~\ref{prop:besselcurves}\,(ii). Since $x_{m'}$ is increasing,
this implies $\nu_m(x)>\nu_{m'}(x)$ and thus the ordering of the zeros $\nu_m(x)$. The proof of Proposition~\ref{prop:besselcurves} directly showed that each $\nu_m$ is
continuous and strictly increasing with values in $(-m,\infty)$. It remains to show that $E(\cdot,x)$ has no zeros
in $\C\setminus\R$ and that all its zeros are simple. Let $N\ge N_0(x)$ be an integer with $N+\frac12>\frac{x^2}{4\log2}$, where $N_0(x)$ is
the constant from Lemma~\ref{lem:circle}. The numbers $\nu_1(x)>...>\nu_N(x)$ are larger than $-N$, and they are smaller
than $N+\frac12$ by Lemma~\ref{lem:E}\,(iv). Thus they are $N$ distinct zeros of $E(\cdot,x)$ in the disk $|\nu|<N+\frac12$. By
Lemma~\ref{lem:circle} this disk contains exactly $N$ zeros, counted with multiplicity. Hence
\begin{equation}\label{eq:allzeros}
\{\nu\in\C:|\nu|<N+\tfrac12,\ E(\nu,x)=0\}=\{\nu_1(x),...,\nu_N(x)\},
\end{equation}
and each of these zeros is simple. Since $N$ can be chosen arbitrarily large, the claim follows.

(ii) By Proposition~\ref{prop:besselcurves}\,(i), the zeros of $E(\nu,\cdot)$ in $(0,\infty)$ are the pairwise distinct numbers
$x_m(\nu)$ with $-m<\nu$, and $x_m(\nu)\le x$ is equivalent to $\nu\le\nu_m(x)$ by monotonicity.

(iii) Let $N$ be as in the proof of (i). By \eqref{eq:allzeros}, the sum of the zeros of $E(\cdot,x)$ in the disk $|\nu|<N+\frac12$
is $\sum_{m=1}^{N}\nu_m(x)$, and \eqref{eq:sigmaN} reads
\[
\bigg|\sum_{m=1}^{N}(\nu_m(x)+m)-\frac{x^2}{4}\bigg|\le\frac{C(x)}{N}.
\]
Letting $N\to\infty$ we obtain \eqref{eq:sumrule}.
\end{proof}

\begin{remark} \label{rem:jacobi}
Proposition~\ref{prop:besselcurves} and Lemma~\ref{thm:bessel} can also be proved without the approximation by Laguerre polynomials, Hurwitz's theorem and without Lemma~\ref{lem:circle}. The idea is to realize the zeros of $\nu\mapsto J_\nu(x)$ as the eigenvalues of a self-adjoint operator. This realization is not new: up to a sign and a unitary transformation, the operator $T_x$ below is the Jacobi operator studied in \cite[Example~20]{StampachStovicek}, see also the references given there.

\begin{proof}[Alternative proof of Proposition~\ref{prop:besselcurves}]
Fix $x>0$. On $\ell^2(\N)$ consider the Jacobi operator
\[
(T_xu)_k=\tfrac x2 u_{k-1}-k u_k+\tfrac x2 u_{k+1},\qquad k\in\N,\quad u_0:=0,
\]
with domain $\{u\in\ell^2(\N):(ku_k)_k\in\ell^2(\N)\}$. We have $T_x=T_0+xK$ with $T_0=-\mathrm{diag}(1,2,3,...)$, where $K$ is the symmetric off-diagonal part with entries $\frac12$, so that $\|K\|\le1$. Hence $T_x$ is self-adjoint, bounded from above and has compact resolvent. If $T_xu=\mu u$ with $u\ne0$, then $u_1\neq0$, since otherwise the first row and the recursion force $u=0$. Consequently all eigenvalues of $T_x$ are simple.

\emph{Step 1: Bessel zeros are eigenvalues.} For $\nu\in\C$ put $u(\nu):=(J_{\nu+k}(x))_{k\in\N}$. This sequence decays faster than any exponential, since $J_{\nu+k}(x)\sim(x/2)^{\nu+k}/\Gamma(\nu+k+1)$ as $k\to\infty$, which follows directly from the series \eqref{eq:E}, and it is not identically zero. The recurrence $J_{\mu-1}(x)+J_{\mu+1}(x)=\frac{2\mu}{x}J_\mu(x)$ \cite[(10.6.1)]{DLMF} gives
\begin{equation}\label{eq:jacobi}
(T_x-\nu)\,u(\nu)=-\tfrac x2\,J_\nu(x)\,e_1 .
\end{equation}
Hence every zero $\nu\in\C$ of $\nu\mapsto J_\nu(x)$ is an eigenvalue of $T_x$, and is therefore real. Conversely, let $\mu$ be an eigenvalue with real eigenvector $u$. Since $T_x$ is symmetric, \eqref{eq:jacobi} gives $0=\langle(T_x-\mu)u,u(\mu)\rangle=-\frac x2J_\mu(x)u_1$, and $u_1\neq0$, so $J_\mu(x)=0$. Thus the zeros of $\nu\mapsto J_\nu(x)$ are exactly the eigenvalues $\nu_1(x)>\nu_2(x)>...$ of $T_x$.

\emph{Step 2: Simplicity.} Let $\nu\in\R$. Differentiating \eqref{eq:jacobi} with respect to $\nu$ gives
\begin{align*}
    (T_x-\nu)\,\partial_\nu u(\nu)-u(\nu)=-\tfrac x2\,\partial_\nu J_\nu(x)\,e_1 .
\end{align*}
The derivative $\partial_\nu u(\nu)$ again lies in the domain of $T_x$. Taking the inner product with $u(\nu)$ and using the symmetry of $T_x$ together with \eqref{eq:jacobi} once more, we obtain
\begin{equation}\label{eq:wronski}
\frac x2\Bigl(J_{\nu+1}(x)\,\partial_\nu J_\nu(x)-J_\nu(x)\,\partial_\nu J_{\nu+1}(x)\Bigr)=\sum_{k=1}^\infty J_{\nu+k}(x)^2 ,
\end{equation}
an identity which was already obtained in \cite[Section~4.2]{StampachStovicek}. At a zero $\nu=\nu_m(x)$ the right-hand side is positive. Hence $\partial_\nu J_\nu(x)\neq0$ and $J_{\nu+1}(x)\neq0$, so all zeros are simple.

\emph{Step 3: Monotonicity.} By the min-max principle, $|\nu_m(x)-\nu_m(x')|\le|x-x'|\,\|K\|\le|x-x'|$, so each $\nu_m$ is continuous. By Step 2 and the implicit function theorem, each $\nu_m$ is real-analytic on $(0,\infty)$. Differentiating $J_{\nu_m(x)}(x)=0$, using $\partial_xJ_\nu(x)=\frac\nu xJ_\nu(x)-J_{\nu+1}(x)$ \cite[(10.6.2)]{DLMF} and applying \eqref{eq:wronski}, gives
\begin{equation}\label{eq:numderiv}
\frac{d}{dx}\nu_m(x)=\frac{J_{\nu+1}(x)}{\partial_\nu J_\nu(x)}=\frac x2\cdot\frac{\,J_{\nu+1}(x)^2}{\sum_{k=1}^\infty J_{\nu+k}(x)^2}>0,\qquad\nu=\nu_m(x).
\end{equation}
Hence every $\nu_m$ is strictly increasing on $(0,\infty)$.

\emph{Step 4: Range.} The min-max principle applied to $T_x=T_0+xK$ gives $|\nu_m(x)+m|\le x$. Hence $\nu_m(x)\to-m$ as $x\searrow0$, and by strict monotonicity of $\nu_m(x)$, we get $\nu_m(x)>-m$ for all $x>0$. To show $\nu_m(x)\to\infty$ as $x\to\infty$, one can use the same argument as in Step~4 of the proof of Proposition~\ref{prop:besselcurves} which uses only strict monotonicity and Lemma~\ref{lem:E}\,(iv). Thus each $\nu_m$ is an increasing bijection from $(0,\infty)$ onto $(-m,\infty)$. With $x_m:=\nu_m^{-1}$, part~(i) of Proposition~\ref{prop:besselcurves} is proven. Part~(ii) follows from $\nu_{m'}<\nu_m$ for $m<m'$, since $\nu_{m'}(x_m(\nu))<\nu_m(x_m(\nu))=\nu=\nu_{m'}(x_{m'}(\nu))$.
\end{proof}

\begin{proof}[Alternative proof of Lemma~\ref{thm:bessel}]
Part (i) is contained in Steps~1-4 above, and part (ii) follows from
Proposition~\ref{prop:besselcurves}\,(i) as before. For (iii), let $x>0$. Again by Steps~1 and~2 above, the
eigenvalues of $T_x$ are the simple eigenvalues $\nu_m(x)$, $m\in\N$, with eigenvectors
$u(\nu_m(x))=(J_{\nu_m(x)+k}(x))_{k\in\N}$. Since $T_x$ is self-adjoint with compact resolvent, the
normalized eigenvectors
\[
\hat u_m(x):=\frac{u(\nu_m(x))}{\|u(\nu_m(x))\|},\qquad m\in\N,
\]
form an orthonormal basis of $\ell^2(\N)$. In terms of these, \eqref{eq:numderiv} reads
\[
\frac{d}{dx}\nu_m(x)=\frac x2\cdot\frac{J_{\nu+1}(x)^2}{\sum_{k=1}^\infty J_{\nu+k}(x)^2}
=\frac x2\,|\langle e_1,\hat u_m(x)\rangle|^2,\qquad\nu=\nu_m(x),
\]
and Parseval's identity gives
\begin{equation}\label{eq:sumderiv}
\sum_{m=1}^\infty\frac{d}{dx}\nu_m(x)=\sum_{m=1}^\infty\frac x2\,|\langle e_1,\hat u_m(x)\rangle|^2=\frac x2\,\|e_1\|^2=\frac x2,\qquad x>0 .
\end{equation}
Each $\nu_m$ is continuously differentiable on $(0,\infty)$ with $\frac{d}{dx}\nu_m>0$ by
Step~3, and $\nu_m(x)\to-m$ as $x\searrow0$ by Step~4. Hence
$\nu_m(x)+m=\int_0^x \frac{d}{dt}\nu_m(t)\,\dd t$ for every $m\in\N$. Since all integrands are positive, the
monotone convergence theorem and \eqref{eq:sumderiv} yield
\[
\sum_{m=1}^\infty (\nu_m(x)+m)=\int_0^x\sum_{m=1}^\infty\frac{d}{dt}\nu_m(t)\,\dd t
=\int_0^x\frac t2\,\dd t=\frac{x^2}{4},
\]
which is \eqref{eq:sumrule}.
\end{proof}

In this approach neither the limit \eqref{eq:limit} nor Lemma~\ref{lem:circle} is needed for the
proof. The limit \eqref{eq:limit} remains as a conceptual link between
Theorems~\ref{thm:landau} and~\ref{thm:weyl}.
\end{remark}

\section{The exact Weyl formula and P\'olya's inequality}\label{sec:weyl}

\begin{proof}[Proof of Theorem~\ref{thm:weyl}]
Let $x=R\sqrt\lambda$. By separation of variables in polar coordinates, the eigenvalues of $H_0^{D_R}$, counted with multiplicity,
are the numbers $j_{|l|,k}^2/R^2$ with $l\in\Z$ and $k\in\N$, where $j_{|l|,1}<j_{|l|,2}<...$ denote the positive zeros of
$J_{|l|}$. By Lemma~\ref{lem:E}\,(v) and \eqref{eq:EJ}, the functions $J_{|l|}$ and $J_l$ have the same positive zeros. Hence
\begin{equation}\label{eq:N0all}
N(\lambda,H_0^{D_R})=\sum_{l\in\Z}\#\{t\in(0,x]:J_l(t)=0\}.
\end{equation}
By Lemma~\ref{thm:bessel}\,(ii), the summand equals $\#\{m\in\N:-m<l\le\nu_m(x)\}$. Interchanging the order of counting gives
\[
N(\lambda,H_0^{D_R})=\sum_{m=1}^{\infty}\#\{l\in\Z:-m<l\le\nu_m(x)\}
=\sum_{m=1}^{\infty} (\lfloor\nu_m(x)\rfloor+m ),
\]
where we used $\nu_m(x)>-m$, see Lemma~\ref{thm:bessel}\,(i). Since the series \eqref{eq:sumrule} converges and its terms are
positive, we have $0<\nu_m(x)+m<1$ and hence $\lfloor\nu_m(x)\rfloor+m=0$ for all sufficiently large $m$. This proves the first
identity in \eqref{eq:weyl}. The second one follows from the sum rule \eqref{eq:sumrule} of Lemma~\ref{thm:bessel}\,(iii). Indeed,
\begin{align*}
\sum_{m=1}^{\infty}(\lfloor\nu_m(x)\rfloor+m )
&=\sum_{m=1}^{\infty}(\nu_m(x)+m )-\sum_{m=1}^{\infty}\{\nu_m(x)\}\\
&=\frac{x^2}{4}-\sum_{m=1}^{\infty}\{\nu_m(x)\}
=\frac{|D_R|}{4\pi}\lambda-\sum_{m=1}^{\infty}\bigl\{\nu_m(R\sqrt\lambda)\bigr\},
\end{align*}
where the last series converges because $0\le\{\nu_m(x)\}\le\nu_m(x)+m$ for any $m\in\N$. 
\end{proof}

\begin{proof}[Proof of Corollary~\ref{cor:polya}]
Let $\lambda>0$ and $x=R\sqrt\lambda$. For all sufficiently large $m$ we have $0<\nu_m(x)+m<1$ and therefore
$\{\nu_m(x)\}=\nu_m(x)+m>0$. Hence the last sum in \eqref{eq:weyl} is strictly positive, which gives the strict P{\'o}lya inequality.  
\end{proof}

\begin{remark}
By \eqref{eq:limit}, the formulas \eqref{eq:weyl} and  \eqref{eq:sumrule} in Lemma~\ref{thm:bessel}\,(iii) are the limits of \eqref{eq:landau} and of Vieta's formula \eqref{eq:vieta_alpha} in Lemma~\ref{cor:roots}\,(iii), respectively, along $\phi=x^2/(4j)$, $j\to\infty$. We stress that the proofs above do not use this limit for the counting function itself. Only the reality and the monotonicity of the zeros are transferred from the Laguerre polynomials to the Bessel function. The alternative proofs of Proposition~\ref{prop:besselcurves} and Lemma~\ref{thm:bessel} avoid the limit \eqref{eq:limit} entirely and decouple the non-magnetic proof from the magnetic proof. 
\end{remark}

\section{Extensions and other consequences} \label{sec:extensions}

In this section we collect some extensions of the results above. Sections~\ref{sec:ab} and~\ref{sec:balls} contain exact counting formulas for two related spectral problems and in Section~\ref{sec:neumann} we formulate analogous formulas for Neumann boundary conditions. Finally, in Section~\ref{sec:flps} we compare the exact Weyl formula \eqref{eq:weyl} with the lattice point count of \cite{FLPS}. 

\subsection{Aharonov-Bohm field}\label{sec:ab}

Fix $\phi\in \mathbb{R}$. Let $H_{AB}^{D_R}=(-i\nabla+\phi A)^2$ with $A(x)=(-x_2,x_1)/|x|^2$ be the magnetic Dirichlet Laplacian on $D_R$ with Aharonov-Bohm field of flux $\phi$ centred at the origin, $\lambda_n^{AB}(D_R)$ its eigenvalues and 
\begin{align*}
        N(\lambda, H_{AB}^{D_R}) &= \#\{n\in\N: \lambda_n^{AB}(D_R) \leq \lambda\}
    \end{align*} 
its counting function. The eigenvalues of $H_{AB}^{D_R}$ are given by $j_{|l-\phi|,k} ^2/R^2$ with $l\in \Z$ and $k\in \N$. Therefore, the following counting formula holds.

\begin{theorem}\label{thm:counting_ab}
Let $R>0$, $\lambda>0$ and $\phi\in \mathbb{R}$. Then
\begin{equation}\label{eq:weyl_ab}
N(\lambda,H_{AB}^{D_R})=\sum_{m: \,\nu_m(R\sqrt{\lambda}) \geq 0}\Bigl(\bigl\lfloor\nu_m(R\sqrt\lambda)+\phi\bigr\rfloor + \bigl\lfloor\nu_m(R\sqrt\lambda)-\phi\bigr\rfloor +1\Bigr).
\end{equation}
\end{theorem}

\begin{proof}
Let $x=R\sqrt{\lambda}$. By Lemma~\ref{thm:bessel}\,(ii) and interchange of the order of counting,
\begin{align*}
N(\lambda,H_{AB}^{D_R})&=\sum_{l\in\Z}\#\{t\in(0,x]:J_{|l-\phi|}(t)=0\}\\
&=\sum_{l\in \Z}\#\{m\in\N: |l-\phi|\le\nu_m(x)\}\\
&=\sum_{m=1}^{\infty}\#\{l\in\Z: |l-\phi|\le\nu_m(x)\}\\
&=\sum_{m: \,\nu_m(x) \geq 0}\Bigl(\bigl\lfloor\nu_m(x)+\phi\bigr\rfloor + \bigl\lfloor\nu_m(x)-\phi\bigr\rfloor +1\Bigr).
\end{align*}
\end{proof}

\subsection{Higher dimensions}\label{sec:balls}

Let $B_R^d = \{ x \in\R^d : |x|< R \}$, $d\geq 2$, be the $d$-dimensional ball of radius $R$, $-\Delta_{B_R^d}^D$ the Dirichlet Laplacian on $B_R^d$, $\lambda_n(B_R^d)$ its eigenvalues and let
    \begin{align*}
        N(\lambda, -\Delta_{B_R^d}^D) &= \#\{n\in\N:\lambda_n(B_R^d) \leq \lambda\}
    \end{align*}
be the associated counting function. The eigenvalues of $-\Delta_{B_R^d}^D$ are given by $j_{l+\frac{d-2}{2},k}^2/R^2$ with $l\in\N_0$ and $k\in\N$, each with multiplicity
    \begin{align*}
        \kappa_{d,l} = \binom{l+d-1}{d-1} - \binom{l+d-3}{d-1},
    \end{align*}
where $\binom nk=0$ for $n<k$. Therefore, the following counting formula holds.

\begin{theorem}\label{thm:counting_d3}
Let $d\geq 2$, $R>0$ and $\lambda>0$. Then
\begin{equation}\label{eq:counting_d}
N(\lambda, -\Delta_{B_R^d}^D)=\sum_{m: \,\nu_m(R\sqrt{\lambda}) \geq \frac{d-2}{2}}A_d\bigl(\bigl\lfloor\nu_m(R\sqrt\lambda)- \tfrac{d-2}{2}\bigr\rfloor \bigr),
\end{equation}
where
\begin{align*}
    A_d(L) = \sum_{l=0}^L \kappa_{d,l} = \binom{L+d-1}{d-1} + \binom{L+d-2}{d-1}.
\end{align*}
\end{theorem}

\begin{proof}
Let $x=R\sqrt{\lambda}$. By Lemma~\ref{thm:bessel}\,(ii) and interchange of the order of counting,
\begin{align*}
N(\lambda, -\Delta_{B_R^d}^D)&=\sum_{l=0}^\infty\kappa_{d,l}\#\{t\in(0,x]:J_{l+ \frac{d-2}{2}}(t)=0\}\\
&=\sum_{l =0}^{\infty} \kappa_{d,l}\#\{m\in\N: l+ \tfrac{d-2}{2}\le\nu_m(x)\}\\
&=\sum_{m: \,\nu_m(R\sqrt{\lambda}) \geq \frac{d-2}{2}} \sum_{l=0}^{\lfloor\nu_m(R\sqrt\lambda)- \tfrac{d-2}{2}\rfloor} \kappa_{d,l}\\
&=\sum_{m: \,\nu_m(R\sqrt{\lambda}) \geq \frac{d-2}{2}}A_d\bigl(\bigl\lfloor\nu_m(R\sqrt\lambda)- \tfrac{d-2}{2}\bigr\rfloor \bigr).
\end{align*}
\end{proof}

\subsection{Neumann boundary conditions}\label{sec:neumann}

We expect that Theorems~\ref{thm:weyl} and~\ref{thm:landau} have natural counterparts for Neumann boundary conditions, which we state in this
subsection. We do not prove these formulas here, but plan to do so in a follow-up article. Numerical evidence substantiates the formulas claimed below.

For $B\ge0$ let $\tilde{H}_B^{D_R}$ denote the Neumann realization of the magnetic Laplacian with constant magnetic field of
strength $B$ on $D_R$, let $\mu_1(D_R,B)\le\mu_2(D_R,B)\le...$ denote its eigenvalues, repeated according to multiplicity, and let
\[
N_<(\mu, \tilde{H}_B^{D_R}) := \#\{n\in\N:\mu_n(D_R,B)<\mu\},\qquad\mu\in\R,
\]
be the left-continuous counting function. Moreover, let $\lceil t\rceil=\min\{k\in\Z:k\ge t\}$ for $t\in\R$, so that
$\lceil t\rceil=t+\{-t\}$.

We begin with the non-magnetic case, where $\tilde{H}_0^{D_R}$ is the Neumann Laplacian on $D_R$. The role of the zeros of
$\nu\mapsto J_\nu(x)$ is taken by the zeros of $\nu\mapsto J_\nu'(x)$, where $J_\nu'$ denotes the derivative of $J_\nu$ with
respect to its argument $x$.

\begin{conjecture}[Neumann Weyl formula for the disk]\label{thm:weyl_neumann}
For every $x>0$ the function $\nu\mapsto J_\nu'(x)$ has infinitely many zeros. All of them are real and simple, and if they are
denoted by $\nu_1'(x)>\nu_2'(x)>...$, then $\nu_m'(x)>-m$ for all $m\in\N$. Moreover, for $R>0$ and $\mu>0$,
\begin{equation}\label{eq:weyl_neumann}
N_{<}(\mu,\tilde{H}_0^{D_R})=\sum_{m=1}^{\infty}\Bigl(\bigl\lceil\nu_m'(R\sqrt\mu)\bigr\rceil+m-1\Bigr)
=\frac{|D_R|}{4\pi}\mu+\sum_{m=1}^{\infty}  \bigl\{-\nu_m'(R\sqrt\mu) \bigr\},
\end{equation}
where only finitely many terms of the first series are different from zero and the second series converges.
\end{conjecture}

Here the prime in $\nu_m'$ is part of the notation and does not indicate a derivative of the function $\nu_m$. Given the first
identity in \eqref{eq:weyl_neumann}, the second one would be equivalent to the sum rule
\begin{equation}\label{eq:sumrule_neumann}
\sum_{m=1}^{\infty}\bigl(\nu_m'(x)+m-1\bigr)=\frac{x^2}{4},\qquad x>0,
\end{equation}
which is the counterpart of \eqref{eq:sumrule}.

As in the Dirichlet case, the remainder in \eqref{eq:weyl_neumann} has a sign. Conjecture~\ref{thm:weyl_neumann} would therefore imply the Neumann side of P\'olya's conjecture for the disk in its strict form, i.e.~for every $R>0$,
\[
N_<(\mu,\tilde{H}_0^{D_R})>\frac{|D_R|}{4\pi}\mu,\qquad\mu>0,
\]
or equivalently,
\[
\mu_{n+1}(D_R,0)<\frac{4\pi n}{|D_R|}, \qquad n\in\N.
\]
The inequalities are strict because the last series in \eqref{eq:weyl_neumann} is positive: for all sufficiently large $m$
we have $\nu_m'(x)\in(-m,-m+1]$, and $\nu_m'(x)\neq-m+1$, since $J_{-l}'=(-1)^lJ_l'$ and $J_l'$ has no zeros in $(0,l)$ for
$l\in\N$, see \cite[\S10.21(i)]{DLMF}. Here, it is important that the counting function is defined as left-continuous by the strict
inequality $<$ instead of $\leq$ to obtain $\mu_{n+1}(D_R,0)$ instead of $\mu_{n}(D_R,0)$ on the left-hand side. The Neumann
P\'olya inequality for the disk is known from \cite{FLPS}, and the first fully analytic proof was claimed in \cite{Li}. A proof via
\eqref{eq:weyl_neumann} does not need to combine several estimates over overlapping ranges of $\mu$, as the proofs in \cite{FLPS}
and \cite{Li} do.

We now turn to the magnetic case and to the counterpart of Theorem~\ref{thm:landau}. For $j\in\N_0$ and $\beta\in\R$ let
\begin{equation}\label{eq:laguerre_neumann}
Q_j(\beta,z)=2z\, \partial_z L_j^{(\beta)}(z) + (\beta-z)\, L_j^{(\beta)}(z), \qquad z\in \R.
\end{equation}
For fixed $z$, this is a polynomial of degree $j+1$ in $\beta$. It takes the role which $P_j(\alpha,z)$ plays for the
Dirichlet problem. Indeed, the eigenfunctions of the magnetic Laplacian on $D_R$ with angular momentum $l\in\Z$ and eigenvalue
$\lambda$ have the radial part $r^{|l|}e^{-Br^2/4}M\bigl(\frac12(l+|l|+1-\frac\lambda B),|l|+1,\frac{Br^2}{2}\bigr)$, compare
\cite[Section~2.1]{BaurWeidl}, and for $\lambda=(2j+1)B$ and $l\le j$ a direct computation based on \eqref{eq:MP} and
\eqref{eq:negparam} shows that the radial derivative of this function vanishes at $r=R$ if and only if $Q_j(-l,\phi)=0$.

\begin{conjecture}[Neumann count at the Landau levels]\label{thm:landau_neumann}
Let $j\in\N_0$. For every $z>0$ the polynomial $\beta \mapsto Q_j(\beta,z)$ has $j+1$ distinct real zeros
$\beta_1^{(j)}(z)>...>\beta_{j+1}^{(j)}(z)$, and they satisfy $\beta_{m}^{(j)}(z)> -m$ for $m=1, ..., j+1$. Moreover, for
$R>0$, $B>0$ and $\phi=BR^2/2$,
\begin{equation}\label{eq:landau_neumann}
N_{<}((2j+1)B,\tilde{H}_B^{D_R})=\sum_{m=1}^{j+1}\Bigl(\bigl\lceil\beta^{(j)}_m(\phi)\bigr\rceil+m-1\Bigr)
=(j+1)\phi+\sum_{m=1}^{j+1} \bigl\{-\beta^{(j)}_m(\phi)\bigr\}.
\end{equation}
In particular,
\begin{equation}\label{eq:twosided_neumann}
(j+1)\phi \leq N_<((2j+1)B,\tilde{H}_B^{D_R}) < (j+1)(\phi+1) ,
\end{equation}
with equality on the left if and only if all $\beta^{(j)}_m(\phi)$ are integers.
\end{conjecture}

Let us comment on why the proofs of Theorems~\ref{thm:weyl} and~\ref{thm:landau} do not carry over directly to
the Neumann case. The zero sets of $(\nu,x)\mapsto J_{\nu}'(x)$ and
$(\beta,z)\mapsto Q_j(\beta,z)$ can in general not be represented as graphs of strictly increasing functions of $\nu$ and
$\beta$: there are (isolated) points in the zero sets with $\partial_x J_{\nu}'(x)=0$ and $\partial_z Q_j(\beta,z) = 0$,
respectively. In numerical examples, one observes that the first zero curve emanates from the point $(0,0)$ to the right, whereas for $m\ge2$ the
$m$-th zero curve emanates from the point $(-(m-1),0)$ to the left and then makes a turn to the right as $x$ and $z$ increase,
respectively. The zero sets are still collections of non-intersecting graphs, now of the functions $\nu_m'$ and $\beta^{(j)}_m$
of $x$ and $z$, respectively, but these functions need not be monotone, and the counting argument has to be adapted to this
situation.

Finally, we expect that Theorem~\ref{thm:counting_d3} also has a counterpart for the Neumann Laplacian $-\Delta_{B_R^d}^N$ on the ball
$B_R^d$, $d\ge2$, as well, which we state without proof. Let $\mu_n(B_R^d)$ denote its eigenvalues and let
\[
N_<(\mu, -\Delta_{B_R^d}^N) := \#\{n\in\N:\mu_n(B_R^d)<\mu\},\qquad\mu\in\R,
\]
be the left-continuous counting function. The eigenvalues are $0$, which is simple with constant
eigenfunction, and the numbers $t_{l,k}^2/R^2$, where $l\in\N_0$, $k\in\N$ and $t_{l,k}>0$ is the $k$-th zero of
\begin{align*}
    t\mapsto tJ_{l+\frac{d-2}{2}}'(t)-\tfrac{d-2}{2}J_{l+\frac{d-2}{2}}(t),
\end{align*}
each with multiplicity $\kappa_{d,l}$. For $x>0$ we denote by $\tilde\nu_1(x)>\tilde\nu_2(x)>...$ the real zeros of
\begin{align*}
    \nu\mapsto xJ_\nu'(x)-\tfrac{d-2}{2}J_\nu(x).
\end{align*}

\begin{conjecture}\label{thm:counting_neumann}
Let $d\ge2$, $R>0$ and $\mu>0$. Then
\begin{equation}\label{eq:counting_neumann_d}
N_<(\mu, -\Delta_{B_R^d}^N)=\sum_{m: \,\tilde\nu_m(R\sqrt{\mu}) > \frac{d-2}{2}}A_d\bigl(\bigl\lceil\tilde\nu_m(R\sqrt\mu)- \tfrac{d-2}{2}\bigr\rceil -1\bigr),
\end{equation}
with $A_d$ from Theorem~\ref{thm:counting_d3}.
\end{conjecture}

\subsection{Relation to lattice point counts of Filonov, Levitin, Polterovich and Sher}\label{sec:flps}

We discuss how the exact Weyl formula \eqref{eq:weyl} relates to the lattice point proof of P{\'o}lya's conjecture for the disk by Filonov, Levitin, Polterovich and Sher in \cite{FLPS}.

For $x>0$ let $G_x(z)=\frac1\pi\bigl(\sqrt{x^2-z^2}-z\arccos\frac zx\bigr)$ for $0\le z\le x$ and $G_x(z)=0$ for $z>x$. Let $\Xi_x\colon[0,\frac x\pi]\to[0,x]$ be the inverse of the strictly decreasing function $G_x|_{[0,x]}$. The proof of P\'olya's conjecture for the disk in \cite{FLPS} compares the counting function with a shifted lattice point count by showing that \cite[Theorem~2.3]{FLPS}
\begin{align} \label{eq:NL_FLPS}
    N(\lambda,H_0^{D_R})\le L(R\sqrt\lambda), \qquad L(x):=\sum_{l\in\Z}\bigl\lfloor G_x(|l|)+\tfrac14\bigr\rfloor .
\end{align}
The key input is \cite[Proposition~3.1]{FLPS} which states that for every real $\nu\ge0$ and $x>0$,
\begin{align} \label{eq:G_upperbound_FLPS}
    \#\{m\in\N:j_{\nu,m}\le x\}\le\bigl\lfloor G_x(\nu)+\tfrac14\bigr\rfloor.
\end{align}
For $m\in\N$ with $m-\frac14\le\frac x\pi$ we have $\lfloor G_x(|l|)+\frac14\rfloor\ge m$ if and only if $|l|\le\Xi_x(m-\frac14)$. Hence
\begin{align} \label{eq:L_rows}
 L(x)=\sum_{m\in\N,\ m-\frac14\le\frac x\pi}\bigl(2\bigl\lfloor\Xi_x(m-\tfrac14)\bigr\rfloor+1\bigr).
\end{align}
On the other hand,
\begin{align} \label{eq:N_rows}
N(\lambda,H_0^{D_R})=\sum_{m:\,\nu_m(x)\ge0}\bigl(2\lfloor\nu_m(x)\rfloor+1\bigr),\qquad x=R\sqrt\lambda ,
\end{align}
for example by Theorem~\ref{thm:counting_ab} with $\phi=0$ or by Theorem~\ref{thm:counting_d3} in the case $d=2$. Thus the lattice count of \cite{FLPS} is obtained by replacing the exact zeros $\nu_m(x)$ by the semiclassical lengths $\Xi_x(m-\frac14)$.

In this picture, the inequality \eqref{eq:G_upperbound_FLPS} becomes an inequality between lengths. By Lemma~\ref{thm:bessel}\,(ii), the left-hand side of \eqref{eq:G_upperbound_FLPS} equals $\#\{m:\nu_m(x)\ge\nu\}$. The right-hand side equals $\#\{m:m-\frac14\le\frac x\pi,\ \Xi_x(m-\frac14)\ge\nu\}$. Both counts are non-increasing in $\nu$. Therefore \eqref{eq:G_upperbound_FLPS} holds for all $\nu\ge0$ if and only if every $m\in\N$ with $\nu_m(x)\ge0$ satisfies $m-\frac14\le\frac x\pi$ and
\begin{align}
   \nu_m(x)\le\Xi_x(m-\tfrac14).  \label{eq:nu_comparison_FLPS}
\end{align}
Since $x_m(\nu)=j_{\nu,m}$ for $\nu>-1$, the inequality \eqref{eq:nu_comparison_FLPS} is exactly the statement $a_{\nu,m}\le j_{\nu,m}$ of \cite[Lemma~3.4]{FLPS}, read along the zero curves.

Subtracting the representations \eqref{eq:L_rows} and  \eqref{eq:N_rows} from each other yields the exact identity
\begin{align*}
    N(\lambda,H_0^{D_R})&=L(x)-E(x),\\
E(x)&=\sum_{m\le n_0(x)}2\Bigl(\bigl\lfloor\Xi_x(m-\tfrac14)\bigr\rfloor-\lfloor\nu_m(x)\rfloor\Bigr)+\sum_{\substack{m>n_0(x)\\ m-\frac14\le\frac x\pi}}\Bigl(2\bigl\lfloor\Xi_x(m-\tfrac14)\bigr\rfloor+1\Bigr),
\end{align*}
where $n_0(x)=\#\{m:\nu_m(x)\ge0\}$. By \eqref{eq:nu_comparison_FLPS}, we have $E(x)\ge0$. This recovers the bound \eqref{eq:NL_FLPS} of \cite{FLPS}, together with an explicit formula for the defect. To conclude P{\'o}lya's inequality for the disk, the authors of \cite{FLPS} show $L(x)<\frac{x^2}4$ by a lattice point argument \cite[Theorem~4.1]{FLPS}. Theorem~\ref{thm:weyl} instead gives $N(\lambda,H_0^{D_R})<\frac{x^2}{4}$ directly.

The same interpretation applies to the other counting problems in \cite{FLPS}. The weighted lattice count for the Dirichlet problem on the ball $B_R^d$ equals $\sum_m A_d\bigl(\lfloor\Xi_x(m-\frac14)-\frac{d-2}{2}\rfloor\bigr)$, summed over $m$ with $\Xi_x(m-\frac14)\ge\frac{d-2}2$. Theorem~\ref{thm:counting_d3} is the same expression with $\Xi_x(m-\frac14)$ replaced by $\nu_m(x)$. For the Neumann problem on the disk, the bound \cite[Proposition~3.1, (3.3)]{FLPS} would be equivalent in the same way to $\nu'_m(x)\ge\Xi_x(m-\frac34)$ for all $m\in\N$ with $m-\frac34\le\frac x\pi$, where $\nu'_m(x)$ are the zeros of $\nu\mapsto J'_\nu(x)$ from Conjecture~\ref{thm:weyl_neumann}.

\section*{Acknowledgements}

The author would like to thank Timo Weidl for helpful comments and remarks.

\section*{AI usage declaration}

The author acknowledges the use of Claude Fable 5.2 for assistance in writing and mathematical reasoning. The underlying idea of Theorem~\ref{thm:landau}, the connection of eigenvalue counts below Landau levels to zeros of Laguerre polynomials, was developed by the author during the work on \cite{BaurWeidl}. An attempted proof of \eqref{eq:twosided} by the author already used the change of perspective of counting zeros of Laguerre polynomials in $\alpha$ instead of $z$ and also identified the Vieta formula \eqref{eq:vieta_alpha} as the relation to obtain the semiclassical phase space volume $j\phi$ in \eqref{eq:landau}. One obstacle that remained unsolved was how to count the eigenvalues of angular momentum $l>0$, which Claude solved by using the relation in Lemma~\ref{lem:Pj}\,(v) and maximal extensions of the zero curves $\zeta_m(\alpha)$ into negative parameter ranges. The idea that one could obtain P{\'o}lya's inequality from \eqref{eq:twosided} by letting $B\to 0$, in particular via Corollary~\ref{cor:flux}\,(i), was also already drafted by the author. The author prompted Claude with the task to complete the proof of \eqref{eq:twosided} given access to unpublished notes by the author. This led to Theorem~\ref{thm:landau}. Theorem~\ref{thm:weyl} was subsequently developed with and proven by the AI, including the sum rule \eqref{eq:sumrule}. All AI-assisted arguments were independently checked by the author, who takes full responsibility for the content of the article.

\end{document}